%% file: ms.tex
\documentclass[3p]{elsarticle}
\pdfoutput=1

\usepackage{hyperref}

\journal{Journal of Computational Science}

\usepackage{natbib}
\usepackage{amsmath}
\usepackage{amsfonts}
\usepackage{graphicx}
\usepackage{amsthm}
\usepackage{moreverb}
\usepackage{amsfonts}
\usepackage{mathrsfs}
\usepackage{float}
\usepackage{fixltx2e}
\usepackage{gensymb}
\usepackage{xfrac}
\usepackage{lmodern}
\usepackage{color}
\usepackage{caption}
\usepackage{subcaption}
\usepackage{epstopdf}
\usepackage{listings}
\usepackage{todonotes}
\usepackage{booktabs}
\usepackage[ruled, lined, linesnumbered, commentsnumbered, longend]{algorithm2e}
\usepackage{subfiles}

\SetKwInOut{KwIn}{Input}
\SetKwInOut{KwOut}{Output}

\newcommand{\R}{\mathbb{R}}
\newcommand{\norm}[1]{\left|\left| #1 \right|\right|}

\begin{document}

\begin{frontmatter}

\title{Reduced Order Modeling for Parameterized Time-Dependent PDEs using Spatially and Memory Aware Deep Learning}

\author{Nikolaj T. M{\"u}cke \corref{mycorrespondingauthor}}
\cortext[mycorrespondingauthor]{Corresponding author}
\ead{nikolaj.mucke@cwi.nl}

\author{Sander M. Bohté}
\ead{S.M.Bohte@cwi.nl}

\author{Cornelis W. Oosterlee}
\ead{C.W.Oosterlee@cwi.nl}

\address{CWI, Science Park 123, 1098 XG Amsterdam, Netherlands}

\begin{abstract}
We present a novel reduced order model (ROM) approach for parameterized time-dependent PDEs based on modern learning. The ROM is suitable for multi-query problems and is nonintrusive. It is divided into two distinct stages: A nonlinear dimensionality reduction stage that handles the spatially distributed degrees of freedom based on convolutional autoencoders, and a parameterized time-stepping stage based on memory aware neural networks (NNs), specifically causal convolutional and long short-term memory NNs. Strategies to ensure generalization and stability are discussed. The methodology is tested on the heat equation, advection equation, and the incompressible Navier-Stokes equations, to show the variety of problems the ROM can handle.
\end{abstract}


\end{frontmatter}


\subfile{Introduction}

\subfile{Parameterized_Time_Dependent_PDEs}

\subfile{Dimensionality_Reduction.tex}

\subfile{Approximating_Paramterized_Time_Evolution_using_Neural_Networks}

\subfile{Results}

\subfile{Conclusion.tex}

\bibliography{references}

\appendix

\subfile{Artificial_Neural_Networks.tex}

\subfile{appendix_results.tex}

\end{document}

%% file: Introduction.tex
\section{Introduction}
Simulations based on first-principles models often form an essential element for understanding, designing, and optimizing problems in, for example, physics, engineering, chemistry, and economics. However, with an increasing complexity of the mathematical models under consideration, it is not always possible to achieve the desired fidelity of such simulations in a satisfactory time frame. This is especially the case when dealing with multi-query and/or real time problems as encountered in uncertainty quantification and model predictive control, where the computational model is typically parameterized. 

There are several approaches to reduce the computation time bottleneck. The arguably most common ones include high-performance computations \cite{hager2010introduction}, high-order discretizations \cite{kopriva2009implementing}, iterative and/or multigrid methods \cite{saad2003iterative, trottenberg2000multigrid}, and reduced order modeling (ROM) \cite{quarteroni2015reduced}. High-performance computing may be costly; the improvements due to high-order discretization strongly depend on the smoothness of solutions at hand, and iterative methods are highly dependent on being able to identify suitable preconditioners. Furthermore, these approaches may suffer from the curse of dimensionality. ROM, a relatively recent research area, is an interesting alternative to the other approaches. 

The ROM solution process is generally divided into two  distinct stages \cite{quarteroni2015reduced}: A so-called ''offline stage'', in which the reduced model is derived, and an ''online stage'', where the reduced model is utilized and solved. Popular choices for the two stages are the proper orthogonal decomposition (POD) model definition, combined with a (Galerkin) projection procedure in the online stage \cite{quarteroni2015reduced, hesthaven2016certified}. Whereas this combination has shown important successes, it has also been shown that the POD and projection approaches perform worse in certain settings, such as for advection-dominated or nonlinear problems. Furthermore, projection-based methods are intrusive, as they require access to the underlying high-fidelity model. Nowadays, it is a reasonable assumption that an industrial model is not directly accessable, and therefore non-intrusive approaches, i.e. approaches that are only based on a series of snapshots of solutions, are increasingly interesting alternatives. 

Machine learning has recently gained the attention from the scientific computing community due to great successes of artificial intelligence in various settings. Specifically Artificial Neural Networks (ANNs), often simply denoted neural networks (NNs), have shown remarkable results in tasks such as image analysis and speech recognition. Much of the success has been boosted further by the availability of open source software frameworks, such as PyTorch  \cite{paszke2017automatic} and Tensorflow \cite{abadi2016tensorflow}, which has made implementation and training possible without expert knowledge and the availability of computation accelerating hardware, such as GPUs, has made training of very large models feasible. These recent advances have accelerated research in especially deep learning, i.e. multilayered NNs, which was not possible few years ago, resulting in many NN architectures specialized in certain tasks, such as time series forecasting and dimensionality reduction, that are able to interpret data with unprecedented accuracy.

Mathematically, NNs have many interesting properties, like universal function approximation \cite{cybenko1989approximation} and pattern recognition that require deep architectures. For these reasons, NNs have gained traction within the mathematics, numerical analysis, and engineering communities either as a replacement or as a supplement to conventional function approximation methods. For an overview of articles, prospects, and future challenges see e.g. \cite{lee2018basic, baker2019workshop, brunton2019machine}. In this paper, we will combine ROM and machine learning in both the offline and the online stages to showcase the potential of using these technologies on conventional problems from scientific computing.

Important work has already been done on the topic of NN-based ROM, which was typically based on Proper Orthogonal Decomposition (POD) for dimensionality reduction and feedforward neural networks (FFNNs) as a parametric map or as a time-stepping scheme \cite{mucke2019reduced,pawar2019deep,pan2018long}. For NN-based dimensionality reduction, see the work in \cite{lee2020model}, where convolutional autoencedoers (CAEs) are used while the time-stepping is done intrusively using multistep methods on a reduced model derived from a Galerkin projection procedure. In \cite{gonzalez2018deep}, autoencoders are utilized for dimensionality reduction and long-short term memory (LSTM) neural networks are used for non-parameterized time-stepping.

Closest to our work are the methodologies described in \cite{fresca2020comprehensive} and \cite{xu2019multi}. In both papers autoencoders are used for dimensionality reduction, and neural networks are used for parameterized time advancement in the latent space. The training of the different networks is done seperately in both papers however. In \cite{fresca2020comprehensive} the dynamics network takes in the parameters and a desired time for the state to be predicted. Hence, the ROM cannot be given an arbitrary initial state and advance in time from there. In \cite{xu2019multi} the time-stepping is further divided into two networks: a causal convolutional neural network (CCNN) is used for encoding the previous states, and a dense FFNN is used to predict the next state based on the encoded dynamics and parameters. Thus, the memory-aware and the parameter neural networks are trained independently from each other. 

In our work, we present a non-intrusive framework, based on deep learning, for computing parameterized spatio-temporal dynamics. The resulting reduced order model is divided into two distinct stages: Firstly, a dimensionality reduction stage based on convolutional autoencoders (CAEs), and secondly a memory-aware NN stage for parameterized time stepping. This methodology utilizes the effectiveness of CAEs as nonlinear dimensionality reduction techniques for spatially distributed data. To discuss the advantages of using CAEs, we make a comparison to the widely used linear counterpart, POD. Specifically we show that POD is a special case of an autoencoder. Furthermore, we present a flexible neural network structure for time-stepping that takes into consideration previous states as well as parameters. The framework is quite general and allows for various types of neural network architectures, hence allowing the researcher to use state-of-the-art techniques that fit the problem at hand. We present and compare two modern time series forecasting architectures, Long Short-Term Memory (LSTM) networks \cite{hochreiter1997long}, and Causal Convolutional Neural Networks (CCNNs) \cite{oord2016wavenet}. Furthermore, we present and discuss a series of approaches to ensure stability and generalization of the time-stepping netowork.

To the best of our knowledge, there is no other work on deep learning-based ROM that is non-intrusive, uses CAEs for dimensionality reduction, has memory-aware and parameterized time-stepping, compares modern time series encoding architectures, and discusses practical approaches to ensure stability and generalization. The result is a flexible offline-online scheme that works for various physical phenomena and can easily be modified according to the specific problem at hand. This makes the presented approach suitable for multi-query problems as they occur in model predictive control, uncertainty quantification, and solving inverse problems.  

The structure of the present paper is as follows. In Section \ref{Parameterized Time-Dependent PDEs} we present parameterized time-dependent PDEs. In Section \ref{Dimensionality Reduction} we discuss dimensionality reduction. Furthermore, we discuss how convolutional autoencoders are used for nonlinear dimensionality reduction. In Section \ref{Approximating Parameterized Time Evolution using Neural Networks} we present the parameterized memory-aware time-stepping neural network. In Section \ref{Results} we showcase the performance on three test problems:  the linear heat equation, a linear advection equation, and the lid driven cavity incompressible Navier-Stokes problem.

%% file: Parameterized_Time_Dependent_PDEs.tex
\section{Parameterized Time-Dependent PDEs}\label{Parameterized Time-Dependent PDEs}

The model under consideration is of the form
\begin{align}\label{partial_diff}
    \partial_t u(t,x;\mu) = F(t,x,u;\mu), \quad u(0,x;\mu) = u_0(x;\mu), 
\end{align}
where $F$ is a (nonlinear) differential operator, $u:\R\times\R^{d}\times\R^{N_p}\rightarrow\R$ or $u:\R\times\R^{d}\times\R^{N_p}\rightarrow\R^d$, $t\in [0,T]$, and $x\in\R^{d}$. Equation \eqref{partial_diff} is a very general parameterized PDE. $\mu$ is to be considered a vector of parameters on which the solution depends. These parameters could be diffusion rate, Reynolds number, parameterize an initial or boundary condition, etc. For technical reasons, the parameter space $\mathcal{P}$ is chosen to be a compact subspace of $\R^{N_P}$  \cite{quarteroni2015reduced}. 

Spatially discretizing \eqref{partial_diff}, using finite elements, finite volumes, finite differences \cite{quarteroni2009numerical}, gives the following finite-dimensional dynamical system,
\begin{align} \label{general_problem}
\partial_t u_h(t,\mu) = F_h(t,u_h(t,\mu);\mu), \quad u_h(0,\mu) = u^0_h(\mu),
\end{align}
$h$ defines the granularity of the discretization, i.e. grid size, number of elements, etc. We will not go into detail regarding these discretizations and it should be assumed that the discretized system is stable and converges to the exact solution when granularity is refined. $u_h(t,\mu)\in \R^{N_h}$ will be referred to as the high-fidelity or full-order solution. 

The manifold of high-fidelity solutions, parameterized by time and the parameters, is called the spatial discrete solution manifold,
\begin{align}
    M_{h} = \left\{ u_h(t,\mu) \: | \: \mu\in \mathcal{P}, \: t\in [0,T]    \right\} \subset \R^{N_h},
\end{align}
Our goal is to approximate this manifold.

Defining a time discretization, $\{t_0,t_1,\ldots,t_{N_t}\}$, $t_n=n\delta t$, and using a time stepping scheme gives us the time discrete approximation of \eqref{general_problem}:
\begin{align} \label{time_discrete_problem}
    u^{n+1}_h(\mu) = F_{h,\delta t}(u^{n}_h;\mu),
\end{align}
where $u^n(\mu)=u(t_n,\mu)$. We will refer to $u^{n}_h(\mu)$ as the state at time step $n$. Note that the discrete time evolution map is not necessarily restricted to only depend on the last state, but can take in several previous states, as is done in e.g. multistep methods, or it could depend on the current state as in implicit methods. We can now define the time-discrete high-fidelity solution manifold:
\begin{align}
    M_{h,\delta t} = \left\{ u^n_h(\mu) \: | \: \mu\in \mathcal{P}, \: n=0,\ldots N_t    \right\} \subset \R^{N_h}.
\end{align}
The subscripts $h$ and $\delta t$ refer to the chosen spatial and time discretizations respectively. $M_{h,\delta t}$ can be seen as the set of discrete state trajectories parameterized by the set of parameters. 

In general, $N_h$ will be very large, which makes advancing the state with \eqref{time_discrete_problem} for many time steps time consuming. This is especially the case when dealing with high-dimensional domains and multiphysics problems. It is indeed a problem when dealing with multi-query problems such as uncertainty quantification and data assimilation or when real-time solutions are of importance as in real-time control settings.

%% file: Dimensionality_Reduction.tex
\section{Dimensionality Reduction} \label{Dimensionality Reduction}

The fundamental idea of dimensionality reduction is that the minimal number of variables necessary to represent the state, also called the intrinsic dimension, of the dynamical system is low compared to the dimension of the high-fidelity model. However, identifying an optimal low-dimensional representation is, in general, not a trivial task. In this section we will give a brief overview of linear dimensionality reduction, particularly, the well-known proper orthogonal decomposition (POD). Then, from the linear outset, we will describe the more general case of nonlinear dimensionality reduction. 

In general, for both linear and nonlinear dimensionality reduction, we assume that a state, $u_h^n(\mu)\in\R^{N_h}$, can be approximated, 
\begin{align} \label{general_dim_reduction}
    u_h^n(\mu) \approx \Phi(u_h^n) = \Phi_{dec} \circ \Phi_{enc} (u_h^n(\mu)),
\end{align}
where $\Phi_{enc} (u_h)\in \R^{N_l}$ with $N_l \ll N_h$. $\Phi_{enc}$ is referred to as the \textit{encoder} and $\Phi_{dec}$ the \textit{decoder}. The encoder transforms the high-dimensional input to a \textit{latent space} of low dimension and the decoder transforms the latent variable back to the high-fidelity space. The latent space is often denoted the \textit{reduced trial manifold}. The state at time step $n$ in the latent space is denoted $u_l^n(\mu) = \Phi_{enc}(u_h^n(\mu))$, and will be referred to as the latent state. 

Ideally, $\Phi$ reconstructs the input perfectly for any given parameters and time step. However, that is, in general, not possible. The precision of the reconstruction is heavily dependent on the dimension of the latent space, as this determines the amount of compression applied. One computes $\Phi$ by  choosing a latent dimension, $N_l$, and then solving the minimization problem,
\begin{align} \label{reduction_minimization}
    \Phi^* = \text{arg}\min_{\Phi} \sqrt{\int_{\mu\in P} \left[\sum_{n=0}^{N_t}\norm{u_h^n(\mu) - \Phi (u_h^n(\mu))}^2_2 \right] d\mu},
\end{align}
where $\norm{\cdot}_2$ denotes the $l^2$-norm.  Theoretically, the reconstruction error should decrease when $N_l$ is increased until the intrinsic dimension of the problem is reached. From thereon, increasing the dimension of the latent space should have very little effect on the reconstruction error. 

There are many ways of solving \eqref{reduction_minimization} \cite{quarteroni2015reduced}. In this paper we focus on a data-driven approach, sometimes referred to as the \textit{method of snapshots}. A snapshot is a high-fidelity solution for a given parameter realization at a certain time. The idea of this approach is to make $N_{train}$ samples from the parameter space and then compute a series of $N_{T}+1$ snapshots, i.e. trajectories, per parameter sample,
\begin{align} \label{M_train}
    M_{N_{train},h,\delta t} = \left\{ u_h^{0}(\mu_1),  \ldots, u_h^{N_t}(\mu_1) ,u_h^{0}(\mu_2),  \ldots, u_h^{N_t}(\mu_2), \ldots , u_h^{0}(\mu_{N_{train}}) , \ldots, u_h^{N_t}(\mu_{N_{train}})   \right\}.
\end{align}
Then \eqref{reduction_minimization} is rewritten into an empirical minimization problem:
\begin{align}\label{reduction_minimization_empirical}
   \Phi^* = \text{arg}\min_{\Phi} \sqrt{\sum_{i=1}^{N_{train}} \sum_{n=0}^{N_t} \left|\left| u_h^{n}(\mu_i) - \Phi \left(u_h^{n}(\mu_i)\right) \right|\right|^2_2}.
\end{align}
The idea is that sampling a finite number of discrete trajectories a sufficient number of times yields a good enough representation of the time-discrete high-fidelity solution manifold. It should be noted that computing \eqref{M_train} is potentially very expensive and even infeasible in some cases.

When a reduction scheme is computed, one can then compute the parameterized trajectories in the latent space, by
\begin{align} \label{time_discrete_problem_reduced}
    u^{n+1}_l(\mu) = F_{l,\delta t}(u^{n}_l;\mu), \quad  u^{0}_l(\mu) = \Phi_{enc}\left(u^{0}_h(\mu)\right),
\end{align}
from which the trajectories in the high-fidelity space can be recovered by $u^{n}_h(\mu) = \Phi_{dec}\left(u^{n}_l(\mu)\right)$. $F_{l,\delta t}$ can be derived in many ways and much time and effort have been put into deriving optimal latent dynamics.

\subsection{Linear Dimensionality Reduction}

In linear dimensionality reduction the strategy is to find a reduced linear trial manifold of low dimension. Since the sought manifold is linear it can be written as the column space, $\text{Col}(V)$ of some matrix, $V\in\R^{N_h\times N_l}$. The column space is the space spanned by the columns of the matrix $V$. From the orthogonal projection theorem, it can be shown that the optimal projection onto a latent linear space is given by
\begin{align} \label{PCA}
    u_h \approx VV^T u_h.
\end{align}
Hence, this is a special case of \eqref{general_dim_reduction} where
\begin{align}
    \Phi = VV^T, \quad \Phi_{enc} = V^T, \quad \Phi_{dec} = V.
\end{align}
This simplification reduces \eqref{reduction_minimization_empirical} to 
\begin{align}\label{reduction_minimization_empirical_linear}
   V^* = \text{arg}\min_{V} \sqrt{\sum_{i=1}^{N_{train}} \sum_{n=0}^{N_t} \left|\left| u_h^{n}(\mu_i) - VV^T u_h^{n}(\mu_i) \right|\right|^2_2},
\end{align}
often accompanied by the constraint that the columns of $V$ are orthogonal, $V^T V=0$. It can be shown that \eqref{reduction_minimization_empirical_linear} has an exact solution \cite{quarteroni2015reduced}. By collecting the snapshots in a \textit{snapshot matrix}, 
\begin{align}
    S = \left[ u_h^{0}(\mu_1) \: | \: \ldots \: | \: u_h^{N_t}(\mu_1) \: |  \: \ldots \: |  \: u_h^{0}(\mu_{N_{train}}) \:| \: \ldots \: | \: u_h^{N_t}(\mu_{N_{train}})   \right],
\end{align}
one can show that the optimal $V\in\R^{N_h\times N_l}$ is given by the first $N_l$ \textit{left singular vectors}. The left singular vectors are computed through the singular value decomposition (SVD),
\begin{align}
    S = U \Sigma Z^T,
\end{align}
where $U$ is a matrix whose columns are the left singular vectors, $Z$ is a matrix whose columns are the right singular vectors, and $\Sigma$ is a diagonal matrix with the singular values on the diagonal. $V$ is then chosen to be the first $N_l$ columns of $U$. This method of obtaining $V$ is the \textit{proper orthogonal decomposition} (POD) \cite{quarteroni2015reduced}, also denoted \textit{principal component analysis} (PCA) \cite{friedman2001elements}.

To obtain $F_{l,\delta t}$ a Petrov Galerkin projection is often performed, which yields
\begin{align}
    F_{l,\delta t}(u^{n}_l;\mu) = W^T F_{h,\delta t}(Vu^{n}_l;\mu).
\end{align}
When $W=V$ it is denoted the Galerkin projection. This approach is \textit{intrusive} which means that direct access to the model, $F_{h,\delta t}$, is required. Furthermore, in the online phase a transformation between the latent space and the high-fidelity space must be performed in each time step in order to be able to evaluate $F_{h,\delta t}(Vu^{n}_l;\mu)$, which slows down the computations. Various methods to circumvent that problem, such the discrete empirical interpolation methods \cite{quarteroni2015reduced},  already exists. However, in recent years there are many studies exploring approximating $F_{l,\delta t}$ with neural networks \cite{pawar2019deep, mucke2019reduced}.

While there are many advantages of a linear reduction scheme, such as the explicit solution to \eqref{reduction_minimization_empirical_linear}, there are, indeed, disadvantages as well. A significant problem is the restriction to a linear trial manifold. The optimal trial manifold, i.e. the trial manifold of the intrinsic dimension, is rarely linear. Especially for advection-dominated and nonlinear problems, it has been shown that a linear reduced approximation does not necessarily lead to significant speed-ups.

\subsection{Nonlinear Dimensionality Reduction}
The extension from linear to nonlinear dimensionality reduction comes naturally and addresses several of the drawbacks of linear dimensionality reduction. The fundamental difference is that we remove the constraint that the latent space has to be a linear manifold. Due to this generalization, we cannot write the projection operator as the matrix product $VV^T$ anymore, but instead we must use the general form in \eqref{general_dim_reduction}, where $\Phi_{enc}$ and $\Phi_{dec}$ can be any type of nonlinear functions. This gives rise to a major difference in solving \eqref{reduction_minimization_empirical}, since no general exact solution exists and therefore \eqref{reduction_minimization_empirical} will be solved numerically.

Even though extra approximation steps have to be introduced in the nonlinear case, the potential gains will, in some cases, outweigh this hurdle. This is due to the fact that with a nonlinear reduction scheme it is theoretically possible to reduce the high-fidelity space down to its intrinsic dimension, $N_P + 1$. However, this relies on the choice of $\Phi$ and the minimization scheme.

A common method for nonlinear dimensionality reduction from the machine learning communities is, among others, kernel PCA. Here, the nonlinear manifold is embedded into a linear space, often of higher dimension, using a predefined nonlinear mapping, $\psi:\R^{N_h}\rightarrow \R^{N_k}$, $N_k>N_h$. From thereon, a linear PCA is performed on the high-dimensional linear data. In order to speed up computations the so-called kernel trick is typically invoked. Utilizing that the nonlinear embedding induces a kernel, $K = k(\psi(x),\psi(y)) = \psi(x)^T\psi(y)$, one can compute the low-dimensional basis without explicitly transforming the data and perform PCA in the high-dimensional space. For more details see \cite{vidal2005generalized}. 

This approach works well in many cases but suffers from one crucial downside: Choosing the nonlinear mapping, $\psi$, or the kernel, $K$, is far from  trivial. There exist no clear guiding principles that work across several cases.

\subsubsection{Autoencoders}
To overcome the problems of other nonlinear dimensionality reduction methods, such as kernel PCA and DEIM, we present Autoencoders (AEs). AEs are a type of NN. For a brief introduction to NNs and the terminology used in this paper, see \ref{appendix:neural_network}. In the context of dimensionality reduction one can interpret an AE as a kernel PCA where the kernel is learned during the training process. Thus, one circumvents the problem of choosing a suitable kernel. Note that this interpretation is merely presented in order to give an intuition of AEs in context of other methods. To further explain the connections between AEs and PCA it is worth noting that a single hidden layer AE with linear activation functions is equivalent to PCA. A single hidden layered AE without bias terms can be written as
\begin{align} \label{one_layer_AE}
    \Phi(u_h^n(\mu);\theta) = T_2\circ T_1(x) = W_2 W_1 u_h^n(\mu), 
\end{align}
where $\theta=\left\{W_1,W_2\right\}$, $T_1:\R^{N_h}\rightarrow \R^{N_l}$, and $T_2:\R^{N_l}\rightarrow \R^{N_h}$ are linear maps and $W_1$ and $W_2$ are matrices. Typically, the mean squared error is chosen as the loss function for AEs, which gives the following minimization problem for the single hidden layer AE:
\begin{align}\label{AE_single_layer_loss}
   \arg \min_{W_1,W_2} \frac{1}{N_{train}N_t}\sum_{i=1}^{N_{train}} \sum_{n=0}^{N_t} \left|\left| u_h^n(\mu_i) - W_2 W_1 u_h^n(\mu_i) \right|\right|^2.
\end{align}
Hence, training a single layer AE is equivalent to solving the PCA minimization problem, eq. \eqref{reduction_minimization_empirical_linear}, without the orthogonality constraint. Conclusively, PCA, eq. \eqref{PCA}, can be considered a special case of an AE. 

By dividing the AE into the encoder and decoder parts and allowing an arbitrary number of layers and nonlinear activation functions, it is easier to understand the similarities to linear dimensionality reduction and why AEs have the potential to perform significantly better. Consider the encoder part with a linear activation in the final layer,
\begin{align}
    \Phi_{enc}(u_h^n(\mu);\theta) =T_L\circ \underbrace{\sigma_{L-1} \circ T_{L-1}  \ldots \circ \sigma_1 \circ T_1 (u_h^n(\mu))}_{\psi_{enc}(u_h^n(\mu);\theta)} = W_L\psi_{enc}(u_h^n(\mu);\theta) = z,
\end{align}
where $\psi_{enc}(x;\theta)=(\psi_{enc}^1(x;\theta),\ldots, \psi_{enc}^{N_e}(x;\theta))\in \R^{N_e}$ and $W_L\in\R^{N_e \times N_l}$. For convenience we ignore bias terms. We see that this corresponds to a nonlinear embedding onto $\R^{N_e}$ and then a projection onto the space spanned by the vectors ${\psi_{enc}^1, \ldots, \psi_{enc}^{N_e}}$. This is quite similar to the idea behind kernel PCA. The difference is that in the AE framework we adjust the nonlinear embedding in the training instead of defining it beforehand. The decoder part is similarly written as
\begin{align}
    \Phi_{dec}(z;\theta) =T_L\circ \underbrace{\sigma_{L-1} \circ T_{L-1}  \ldots \circ \sigma_1 \circ T_1 (z)}_{\psi_{dec}(z)} = W_L\psi_{dec}(z;\theta) = \tilde{u}_h^n(\mu),
\end{align}
where $W_L\in\R^{N_d\times N_h}$.

Note that this is merely a brief discussion of the topic of AEs aiming to give an intuitive understanding. For more details see \cite{yu2019understanding}. 

\subsubsection{Convolutional Autoencoders}

Convolutional autoencoders (CAEs) are a special type of AEs utilizing convolutional layers instead of dense layers. A brief introduction to convolutional neural networks (CNNs) can be found in \ref{appendix:neural_network}. It can be shown that dense and convolutional neural networks are equivalent regarding approximation rates \cite{petersen2020equivalence}, which means that theoretical approximation results for dense NNs translate directly to CNNs. For practical purposes, however, convolutional layers are often to be preferred due to especially the following two properties:
\begin{itemize}
    \item \textit{Shared weights}, which in practice makes the affine transformations very sparse and enables location invariant feature detection.
    \item \textit{Local connections}, which utilizes that spatial nodes close to each other are highly correlated.  
\end{itemize}
An additional advantage is that it is straightforward to handle multiple spatially distributed states. These occur in coupled PDEs such as the Navier Stokes equations where one is dealing with both the $x-$, $y-$ and $z-$components of the velocity field as well as the pressure field. In the framework of CAEs these can all be included by interpreting them as different channels. This enables the possibility of including multiple spatial states without increasing the number of weights in the neural network significantly. The connection between PDEs and CNNs has already been made, see e.g. \cite{ruthotto2019deep}. 

\begin{figure}
    \centering
    \includegraphics[width=0.9\linewidth]{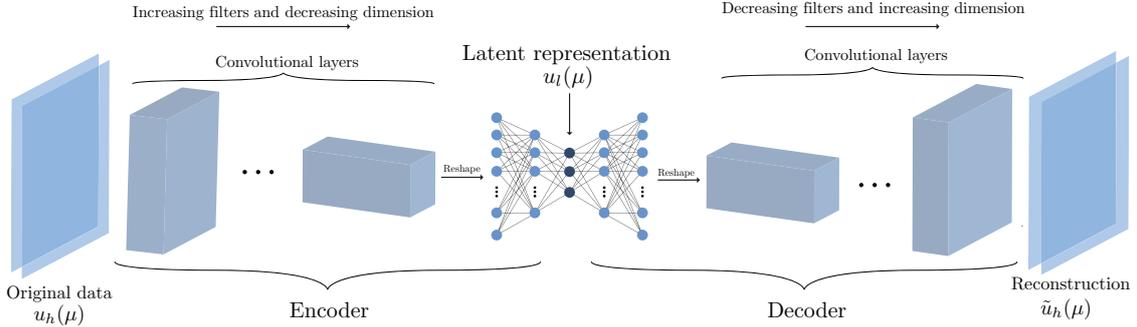}
    \caption{Illustration of a convolutional autoencoder}
    \label{fig:conv_AE}
\end{figure}

In Figure \ref{fig:conv_AE} one sees an illustration of a CAE. The encoding consists of a series of convolutional layers with an increasing number of filters and decreasing dimension, effectively down sampling the number of degrees of freedom, followed by dense layers. Similarly, the decoding consists of series of dense layers followed by a series of deconvolutional layers with a decreasing number of filters and increasing dimension, effectively up sampling. The down sampling is often performed by utilizing pooling layers or strides larger than one. 

It is worth noting that computing the decoder, $\Phi_{dec}$, of a CAE in the training phase is effectively solving an inverse problem. Inverse problems are, in general, ill-posed and therefore require some sort of regularization. $L^2$-regularization, often referred to as weight decay, is frequently used, and results in the following minimization problem to solve:
\begin{align}\label{CAE_loss}
   \arg \min_{\theta} \frac{1}{N_{train}N_t}\sum_{i=1}^{N_{train}} \sum_{n=0}^{N_t} \left|\left| u_h^n(\mu_i) - \Phi(u_h^n(\mu_i);\theta) \right|\right|^2 + \alpha \norm{\theta}_2^2,
\end{align}
where $\alpha$ is a hyperparameter to be tuned. Besides ensuring well-posedness the term also ensures generalization.  

%% file: Approximating_Paramterized_Time_Evolution_using_Neural_Networks.tex
\section{Approximating Parameterized Time Evolution using Neural Networks}\label{Approximating Parameterized Time Evolution using Neural Networks}
In the previous section we presented the general framework for nonlinear dimensionality reduction and showcased how convolutional autoencoders fit into this framework. In this section we aim to explain how neural networks will be utilized for approximating the dynamics in the latent space. For a brief review of the relevant types of neural networks, see \ref{appendix:neural_network}. Neural networks have already shown to be able to approximate dynamical systems \cite{xu2019multi, fresca2020comprehensive, erichson2019physics}. For this reason, it makes sense to test out various architectures for this exact purpose. 

We aim to approximate the dynamics in the latent space non-intrusively by a function, $\Psi \approx F_{l,\delta t}$: 
\begin{align}
    u^{n+1} = \Psi(u^{n}),
\end{align}
where $\Psi$ is a neural network. The approximated latent states will be denoted, $\tilde{u}^n_l(\mu)$, to distinguish from the encoded high-fidelity state, $u_l^n(\mu)=\Phi_{enc}(u_h^n(\mu))$. Thereby, we aim to achieve:
\begin{align}
\begin{split}
    &\left\{ \tilde{u}_l^{0}(\mu_1),  \ldots, \tilde{u}_l^{N_t}(\mu_1) , \ldots , \tilde{u}_l^{0}(\mu_{N_{train}}) , \ldots, \tilde{u}_l^{N_t}(\mu_{N_{train}})   \right\} \\
    \approx &\left\{ u_l^{0}(\mu_1),  \ldots, u_l^{N_t}(\mu_1) , \ldots , u_l^{0}(\mu_{N_{train}}) , \ldots, u_l^{N_t}(\mu_{N_{train}})   \right\}
\end{split}
\end{align}

\subsubsection*{Taking Longer Steps}
Using high-fidelity methods for time-stepping often includes some restrictions on the step size in order for the scheme to be stable. An example is the Courant–Friedrichs–Lewy (CFL) condition for advection-dominated problems \cite{leveque2007finite}. With our strategy, where we aim to learn a neural network representation of the time evolution map, there is no immediate connection between step size and stability. Therefore, in order to speed up online computations, the neural network can be trained to learn to take steps of size $s\delta t$. Hence, $\Psi \approx F_{l,s\delta t}$. 

In the offline phase the high-fidelity trajectories are still computed with step size $\delta t$, to ensure stability, but only every $s$'th step is used for training the NN:
\begin{align}
    \left\{ u_h^{0}(\mu),u_h^{1}(\mu),u_h^{2}(\mu)  \ldots, u_h^{N_t}(\mu)  \right\} \mapsto \left\{ u_h^{0}(\mu), u_h^{s}(\mu), u_h^{2s}(\mu)  \ldots, u_h^{N_t}(\mu)  \right\}.
\end{align}
In general, two states one time-step apart, say $u_h^{n}(\mu)$ and $u_h^{n+1}(\mu)$, are highly correlated. In practice the means that we gain very little extra information by using both in the training of the NN. Therefore, it makes sense only use every $s$'th step to save memory and speed up the training. Hoeever, the number $s$ must be chosen according to various factors, like the requested detail of the dynamics in the online phase. It should further be kept in mind that larger $s$ results in a more complicated map to learn, and thus complicates the training. 

For simplicity we will use the notation $\left\{ u_h^{0}(\mu),u_h^{1}(\mu),u_h^{2}(\mu)  \ldots, u_h^{N_t}(\mu) \right\}$ when referring to the trajectory used for training the neural network. 

\subsubsection*{Approximating the State vs. Residual}
At first glance, it makes sense to train a neural network to approximate $u_l^{n+1}$ directly given $u^{n}$. However, it is shown in \cite{pawar2019deep} and \cite{gin2019deep} that learning the residual instead of the next state often improves the accuracy. Hence, we consider the case
\begin{align}
    u_l^{n+1} = \Psi(u^{n}) = u_l^{n} + R(u_l^{n}),
\end{align}
where $R$ is being approximated by a neural network. This practically makes $\Psi$ what is often referred as a residual neural network. 

\subsubsection*{Incorporating Memory}
In \cite{pawar2019deep} and \cite{gonzalez2018deep} the potential benefits of not only using the present state but also incorporating several previous time steps for the future predictions were shown. Therefore, we now consider
\begin{align}
    u_l^{n+1} =\Psi(u_l^{n},u^{n-1},\ldots,u_l^{n-\xi}) = u_l^{n} + R(u_l^{n},u_l^{n-1},\ldots,u_l^{n-\xi}),
\end{align}
where $\xi$ is the number of previous states included as input into the residual computation by the NN. The idea of incorporating several previous timesteps can loosely be compared to linear multistep methods where the order of approximation can be increased by using several previous steps \cite{leveque2007finite}. In contrast to linear multistep methods, NNs incorporate the previous time steps in a nonlinear fashion.

\begin{figure}
\centering
  \includegraphics[width=.9\linewidth]{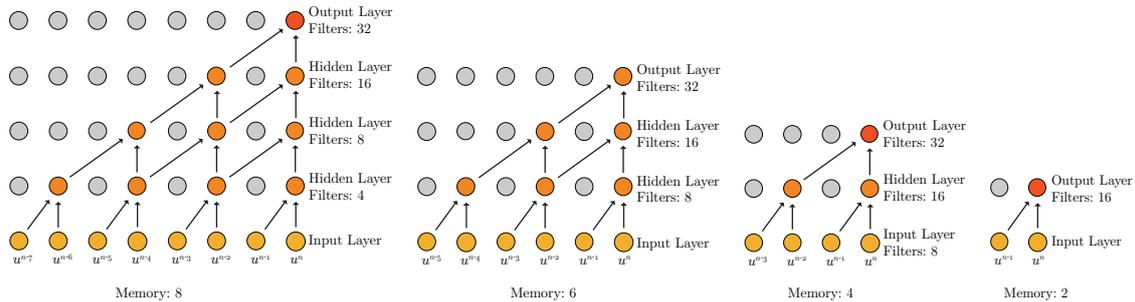}  
\caption{CCNN network architectures for varying memory.}
\label{fig:CCNN_types}
\end{figure}

We consider two different types of networks in this paper: LSTM, and the CCNN. The two types of neural networks take varying computational time to train, have varying numbers of parameters, and vary in regards to how they interpret memory. In the appendix, there is a short description of the two types. Regarding the CCNNs, there are different ways to include memory. In this paper, we have chosen to include memory in shape of adding more layers, see Figure \ref{fig:CCNN_types} for examples for $\xi=8$, $\xi=6$, $\xi=4$, and $\xi=2$.

\subsubsection*{Parameterized Dynamics}
\begin{figure}
    \centering
    \includegraphics[width=0.75\linewidth]{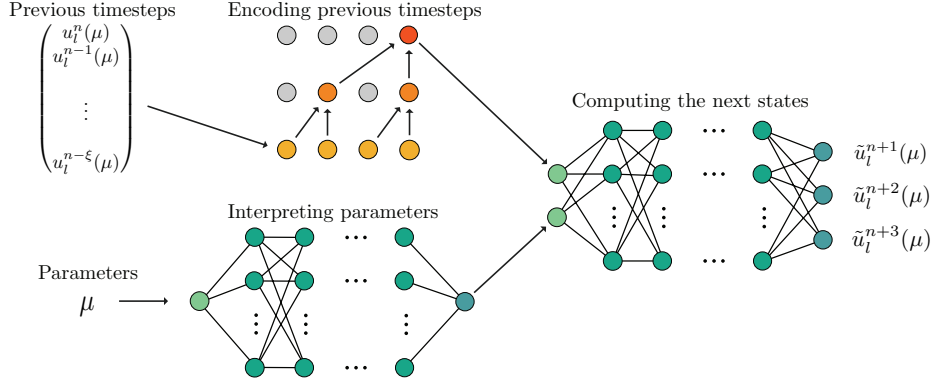}
    \caption{Illustration of the parallel neural network structure. The "Encoding previous timesteps" part is visualized using the CCNN, but it should be noted that an LSTM network (or any suitable time series encoder) could be put in its place. }
    \label{fig:network_parallel}
\end{figure}

We aim to simulate parameterized trajectories of the latent dynamics. Hence, we need to incorporate the parameters as input to the residual computation. For now, we consider constant parameters, but note that it should be possible to incorporate time-dependent parameters. For this reason, the parameters do not need to be part of the memory aware section of the network. We propose a parallel architecture consisting of two branches combining into one: One branch interpreting the last $\xi$ states and one branch processing the parameters. The two branches then connect and provide one final prediction for the residual. Having a single neural network incorporating the previous states and the parameters enables simultaneous training of the two branches. See Figure \ref{fig:network_parallel} for an illustration of the network structure. This ensures that the learned latent features from both branches are optimal with respect to predicting the next state. This is in contrast to what is done in \cite{xu2019multi} where the memory and the parameters are incorporated into two completely separate networks. 

For the parameter branch we simply make use of a dense FFNN. There is no immediate reason to believe that more complicated architectures are necessary, since we are not dealing with time-dependent nor high-dimensional or continuously spatially varying input. Note, however, that there is no reason to believe that this methodology will not work if the FF network in the parameter branch is replaced with a memory aware network in more advanced settings. 

The training of the full time-evolution network is done by minimizing the loss function
\begin{align}
    L(u_l,\mu;\theta) = \frac{1}{N_{train}N_{t}}\sum_{i=1}^{N_{train}}\sum_{n=\xi}^{N_T} \norm{u^{n+1}(\mu_i)-\Psi(u_l^{n}(\mu_i),\ldots,u_l^{n-\xi}(\mu_i),\mu_i;\theta)}^2_2,
\end{align}
with respect to the NN parameters $\theta$. 

\subsubsection*{Imposing Stability and Generalization}
It is well-known that NNs do not necessarily generalize well beyond the training data without some kind of regularization. Combining that with the general risk of having instability in discrete dynamical systems makes it crucial to address these problems during the training. 

The arguably most common technique is to add $L^1$- or $L^2$-regularization to the loss function. Furthermore, specifically for dynamical systems, it has been shown in \cite{erichson2019physics} and \cite{pan2018long} that regularizing the eigenvalues of the Jacobian of the dynamics with respect to the state variable, $D_u \Psi$, does improve long term predictions. In short this is related to linear and Lyapunov stability analysis of dynamical systems, which deals with the analysis of when a system is unstable. Hence, we propose adding the term $\norm{D_u \Psi}_2$, which is the matrix 2-norm, i.e. the spectral radius of the Jacobian of $\Psi$, to the loss function. In practice, by utilizing the relation
\begin{align}
    \norm{D_u \Psi}_2 \leq \norm{D_u \Psi}_F,
\end{align}
we instead add the computationally much cheaper Frobenius norm.

It can empirically be shown that the long term predictions are significantly better if the network takes several steps at a time instead of a single one. Hence, we modify the output of the NN to
\begin{align}
    R\left(u_l^{n}(\mu),\ldots,u_l^{n-\xi}(\mu),\mu;\theta\right) = \left[R_1, R_2, \ldots, R_\zeta\right]^T,
\end{align}
which gives future predictions,
\begin{align}
\begin{bmatrix}
u_l^{n+1}(\mu) \\ \vdots\\  u_l^{n+\zeta}(\mu) 
\end{bmatrix} 
= \begin{bmatrix}
\Psi_1\left(u_l^{n}(\mu),\ldots,u_l^{n-\xi}(\mu),\mu;\theta\right) \\ \vdots\\  \Psi_{\zeta}\left(u_l^{n}(\mu),\ldots,u_l^{n-\xi}(\mu),\mu;\theta\right)
\end{bmatrix} = 
u_l^{n}(\mu) + 
\begin{bmatrix}
R_1\left(u_l^{n}(\mu),\ldots,u_l^{n-\xi}(\mu),\mu;\theta\right) \\ \vdots\\  R_{\zeta}\left(u^{n}(\mu),\ldots,u_l^{n-\xi}(\mu),\mu;\theta\right)
\end{bmatrix} 
\end{align}
Empirically we see that this modification keeps the prediction from exploding for longer time and it reduces spurious oscillations. 

The final loss function for the dynamics NN is given by:
\begin{align}
\begin{split}
    L(u,\mu;\theta) &=\frac{1}{N_{train}N_{t}} \sum_{i=1}^{N_{train}}\sum_{n=\xi}^{N_T} \norm{\sum_{k=1}^\zeta \left[u_l^{n+k}(\mu_i)- \Psi_k\left(u_l^{n}(\mu_i),\ldots,u_l^{n-\xi}(\mu_i),\mu_i;\theta\right)\right]}_2^2\\
    &+ \underbrace{\beta_1 \norm{\theta}_2^2}_{\text{Weight decay}} + \underbrace{\beta_2 \norm{D_u R}_F}_{\text{Jacobian regularization}},
\end{split}
\end{align}

\subsection{The Complete Scheme}
Putting it all together, we have a scheme subdivided into two parts that are trained independently: The CAE, and  the time evolution. The whole process is divided into an online phase and an offline phase. 

\begin{figure}
\centering
\begin{subfigure}{1.\textwidth}
  \centering
  \includegraphics[width=1.\linewidth]{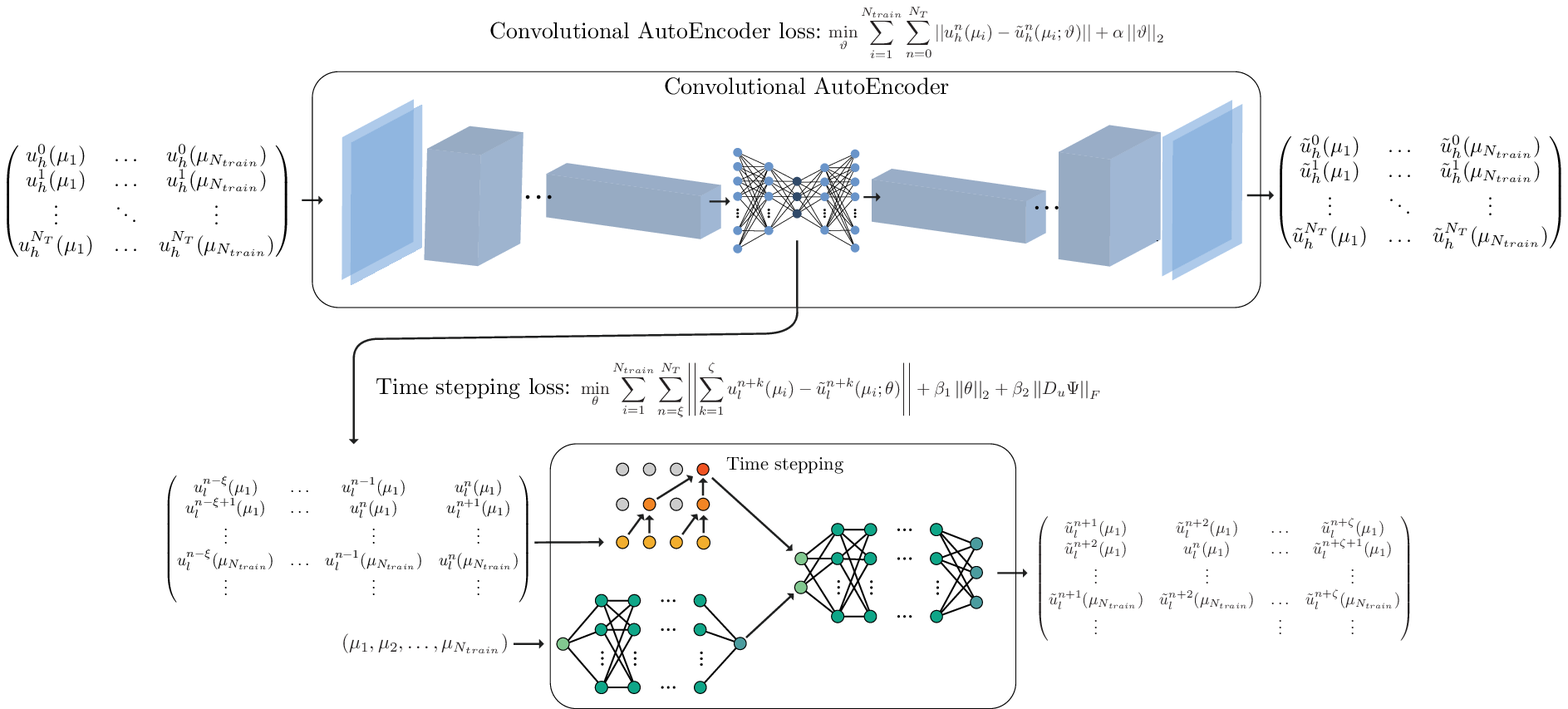}
  \caption{Illustration of the offline stage.}
  \label{fig:offline_stage}
\end{subfigure} \\
\begin{subfigure}{1.\textwidth}
  \centering
  \includegraphics[width=1.\linewidth]{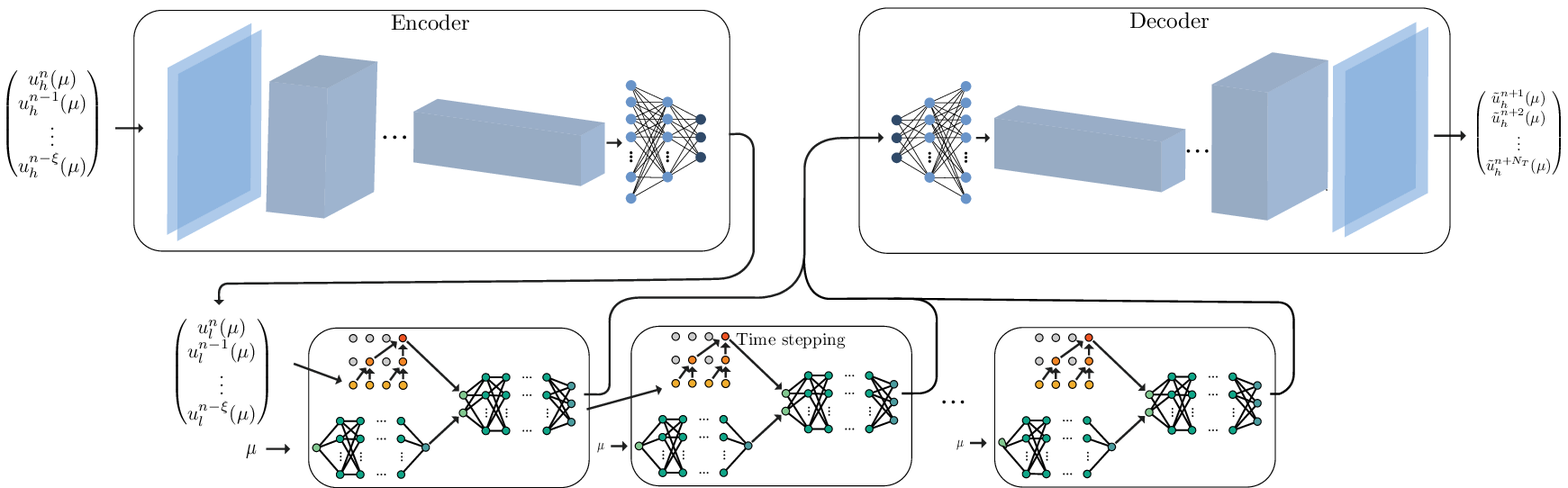}
  \caption{Illustration of the online stage.}
  \label{fig:online_stage}
\end{subfigure}
\caption{}
\label{fig:stages}
\end{figure}

In the offline phase the CAE is trained on a series of high-fidelity snapshots in order to identify a nonlinear reduced trial manifold. Then the CAE is used to reduce the high-fidelity snapshots to the latent space. The latent space trajectories are then used to train the time evolution NN. The training of the two neural networks is visualized in Figure \ref{fig:offline_stage} and outlined in Algorithm \ref{algo:offline_train}. Note that in steps 3 and 5, where the autoencoder and the time evolution network, respectively, are being trained, the considerations mentioned in Appendix \ref{appendix:neural_network} have to be included. Here, we refer to early-stopping, multiple-initialization, choice of optimizer, etc. In Algorithm \ref{algo:offline_tuning} an algorithm to automatically choose the latent dimension, number of training trajectories, memory, and future steps per iteration is presented. Note that this is a very simple approach to tune the network. More advanced methods such as Bayesian optimization or reinforcement learning could be utilized here. Furthermore, it is worth noting that we can, assuming no time constraints, generate as many training samples as necessary. 

In the online phase the first $\xi$ time steps of the state, computed with a high-fidelity scheme for a given parameter realization $\mu$, are projected onto the latent space using the encoder part of the CAE. From there, the time evolution NN computes the parameterized latent space trajectories iteratively. The latent space trajectories are then transformed to the high-fidelity space using the decoder of the CAE. The online stage is visualized in Figure \ref{fig:online_stage} and described in pseudo code in Algorithm \ref{algo:online_train}.

\begin{algorithm} 
\caption{Offline Stage - Training}\label{algo:offline_train}
\KwIn{$N_l$, $\zeta$, $\xi$, $N_{train}$.}
\KwOut{Trained CAE, $\Phi$, and Time Evolution Network, $R$, and test error, E.}

Sample $N_{train}$ parameter samples from the parameter space. \\
Generate high-fidelity trajectories,
$$\left\{ u_h^{0}(\mu_1),  \ldots, u_h^{N_t}(\mu_1) , \ldots , u_h^{0}(\mu_{N_{train}}) , \ldots, u_h^{N_t}(\mu_{N_{train}})   \right\}.$$ \\
Train CAE, $\Phi=\Phi_{dec}\circ\Phi_{enc}$, with latent space dimension $N_l$, by minimizing
\begin{align}
   \arg \min_{\theta} \frac{1}{N_{train}N_t}\sum_{i=1}^{N_{train}} \sum_{n=0}^{N_t} \left|\left| u_n(\mu_i) - \Phi(u_n(\mu_i);\theta) \right|\right|^2.
\end{align} \\
Encode high-fidelity trajectories to get latent state space trajectories
$$\left\{ \Phi_{enc}\left(u_h^{0}(\mu_1)\right),  \ldots, \Phi_{enc}\left(u_h^{N_t}(\mu_1)\right) , \ldots , \Phi_{enc}\left(u_h^{0}(\mu_{N_{train}})\right) , \ldots, \Phi_{enc}\left(u_h^{N_t}(\mu_{N_{train}})\right)   \right\}.$$ \\
Train time evolution network, $R$, to take the last $\zeta$ states and output the residuals for the next $\xi$ states, by minimizing
\begin{align*}
\frac{1}{N_{train}N_{t}} \sum_{i=1}^{N_{train}}\sum_{n=\xi}^{N_T} \norm{\sum_{k=1}^\zeta \left[u_l^{n+k}(\mu_i)- \Psi_k\left(u_l^{n}(\mu_i),\ldots,u_l^{n-\xi}(\mu_i),\mu_i;\theta\right)\right]}_2^2
    + \beta_1 \norm{\theta}_2^2+ \beta_2 \norm{D_u R}_F,
\end{align*} \\
Estimate error on a test set.
\end{algorithm}

\begin{algorithm}
\caption{Offline Stage - Tuning} \label{algo:offline_tuning}
\KwIn{Desired test error, $E^*$}
\KwOut{Optimal latent dimension, $N_l$, $\zeta$, $\xi$, and $N_{train}$.}
Initialize $E = \infty$, $N_l=1$, $\zeta=0$, $\xi=1$, $N_{train}$. \\
\While{$E^*<E$}
{
    Train $\Phi$ and $R$ and compute test error, $E$, using Algorithm \ref{algo:offline_train}. \\
    Update $N_l$, $\zeta$, $\xi$, $N_{train}$ and according to some update rule. 
}
\end{algorithm}

\begin{algorithm}
\caption{Online Stage}  \label{algo:online_train}
\KwIn{$\Phi_{dec}$, $R$, $\mu$, $u^0_h(\mu),\ldots, u^{\xi}_h(\mu)$}
\KwOut{Approximated trajectory in high-fidelity space.}

Encode the initial $\xi$ high-fidelity states,
\begin{align*}
    (u^0_l(\mu),\ldots, u^{\xi}_l(\mu)) = \left(\Phi_{enc}\left(u^0_h(\mu)\right), \ldots, \Phi_{enc}\left(u^\xi_h(\mu)\right)    \right)
\end{align*} \\
Compute approximated latent trajectory by iterating,
\begin{align*}
    \left(\tilde{u}^{n}_l(\mu),\ldots,\tilde{u}^{n+\zeta}_l(\mu)\right) = \tilde{u}^{n}_l(\mu) + R\left(\tilde{u}^{n}(\mu),\ldots,\tilde{u}^{n-\xi}(\mu),\mu;\theta\right),
\end{align*}
until desired end time has been reached. \\
Decode approximated latent space trajectories to high-fidelity space:
\begin{align*}
\left\{ \tilde{u}_h^{0}(\mu),  \ldots, \tilde{u}_h^{N_t}(\mu) \right\} = 
    \left\{ \Phi_{dec}\left(\tilde{u}_l^{0}(\mu)\right),  \ldots, \Phi_{dec}\left(\tilde{u}_l^{N_t}(\mu)\right)  \right\}.
\end{align*}
\end{algorithm}

%% file: Results.tex
\section{Results} \label{Results}

The aim of this section is to showcase how well our frameworks perform for various parameterized PDE problems. Furthermore, we show how the various approaches, regularizations, and parameters affect the performance.  

To assess the performance measure the error on $N_{test}$ test trajectories for parameters values, $\left\{\mu_1,\ldots,\mu_{N_{test}}\right\}$, that the NNs have not seen in the training phase. we measure the mean relative error (MRE) at every time step and take the mean over multiple runs of the test cases:
\begin{align}
    \text{MRE}(u^n_h(\mu_i),\tilde{u}^n_h(\mu_i)) = \frac{1}{N_{test}} \sum_{i=1}^{N_{test}}
    \frac{\norm{u^n_h(\mu_i) - \tilde{u}^n_h(\mu_i)}_2^2}{\norm{u^n_h(\mu_i)}_2^2},
\end{align}
where 
\begin{align}
    \norm{u^n_h(\mu_i)}_2^2 = (u_h^n(\mu_i))^T u_h^n(\mu_i).
\end{align}
Besides showing the MRE, we also show the standard error:
\begin{align}
    \text{Standard Error} = \frac{\sigma}{\sqrt{{N_{test}}}},
\end{align}
where $\sigma$ is the variance of MRE. With this measure we can assess if the trained NN performs similarly on all the test data. I.e. we empirically show robustness and generalization.

For all problems we compare the reconstruction error of the CAE with the reconstruction error of a POD approach. Note that regarding POD the reconstruction error is the same as the projection error. The measure for the POD error is the MRE as for the CAE. 

\subsection*{Neural Network Setup}
All neural networks are implemented in Tensorflow 2.0 \cite{abadi2016tensorflow} in Python. The training is performed in the Google Colab framework on NVIDIA Tesla P100 GPUs. \\ 
The neural network architecture configurations can be found in \ref{conv_ae_appendix} and \ref{evolution_NN_appendix}. \\
\\
\textbf{Remark.} We only present results on dimensionality reduction using convolutional autoencoders and compare them to POD. It should be noted that dense autoencoders were also tested and showed significantly worse results. \\ 
\\
\textbf{Remark.} As mentioned earlier, we only consider FFNNs, LSTMs, and CCNNs for the time -stepping. We also studied other achitectures, such as neural ODEs \cite{chen2018neural}, gated recurrent units (GRUs), and simple recursive neural networks. However, we chose to not include those results. Neural ODEs performed significantly worse and the training took much longer time. GRUs performed similarly to LSTMs and simple recursive neural networks performed slightly worse. 

\subsection{Heat Equation}
\begin{figure}
\centering
\begin{subfigure}{.4\textwidth}
    \centering
    \includegraphics[width=1.\textwidth]{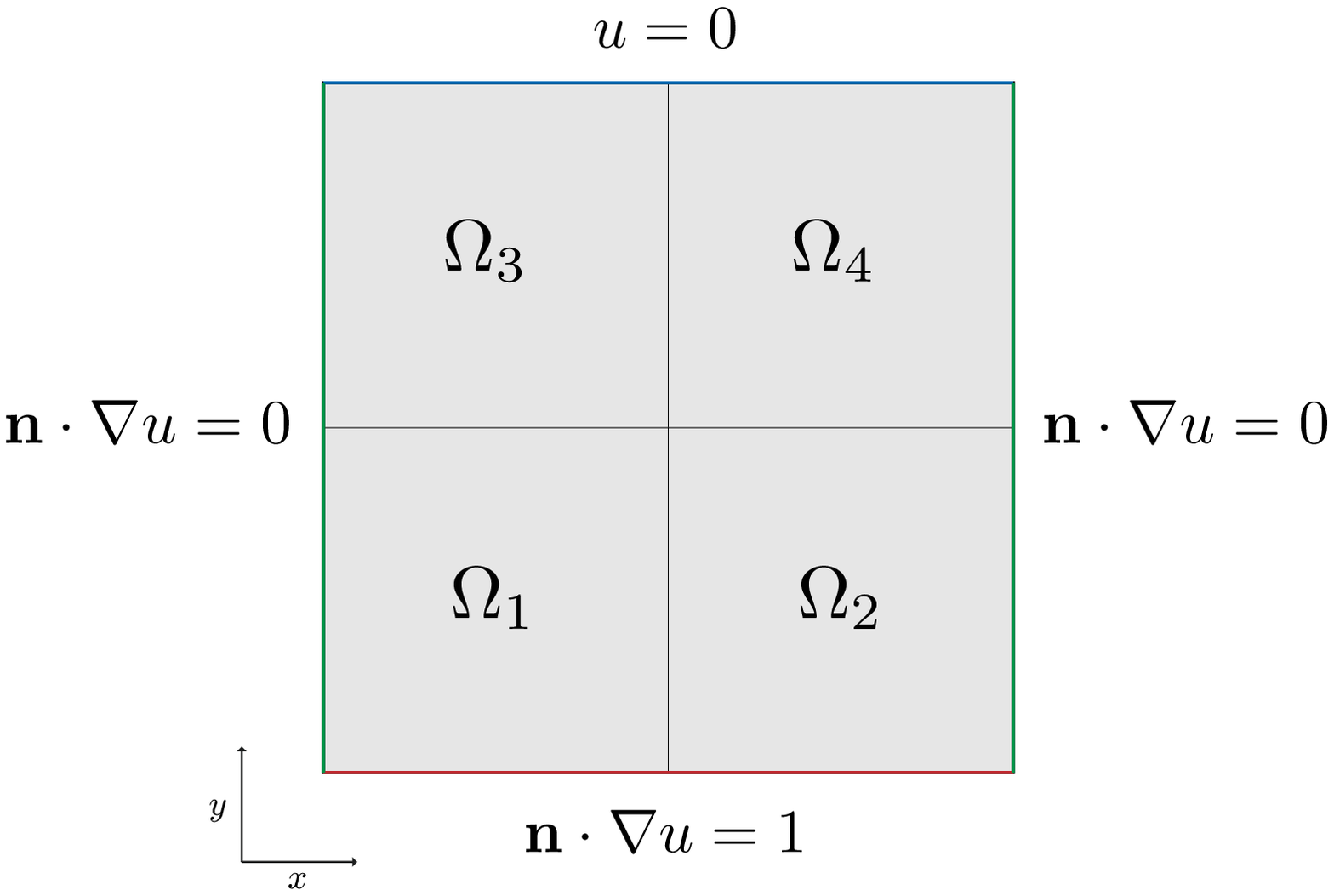}
    \caption{Heat equation setup.}
    \label{fig:heat_setup}
\end{subfigure}
\begin{subfigure}{.4\textwidth}
    \centering
    \includegraphics[width=1.\textwidth]{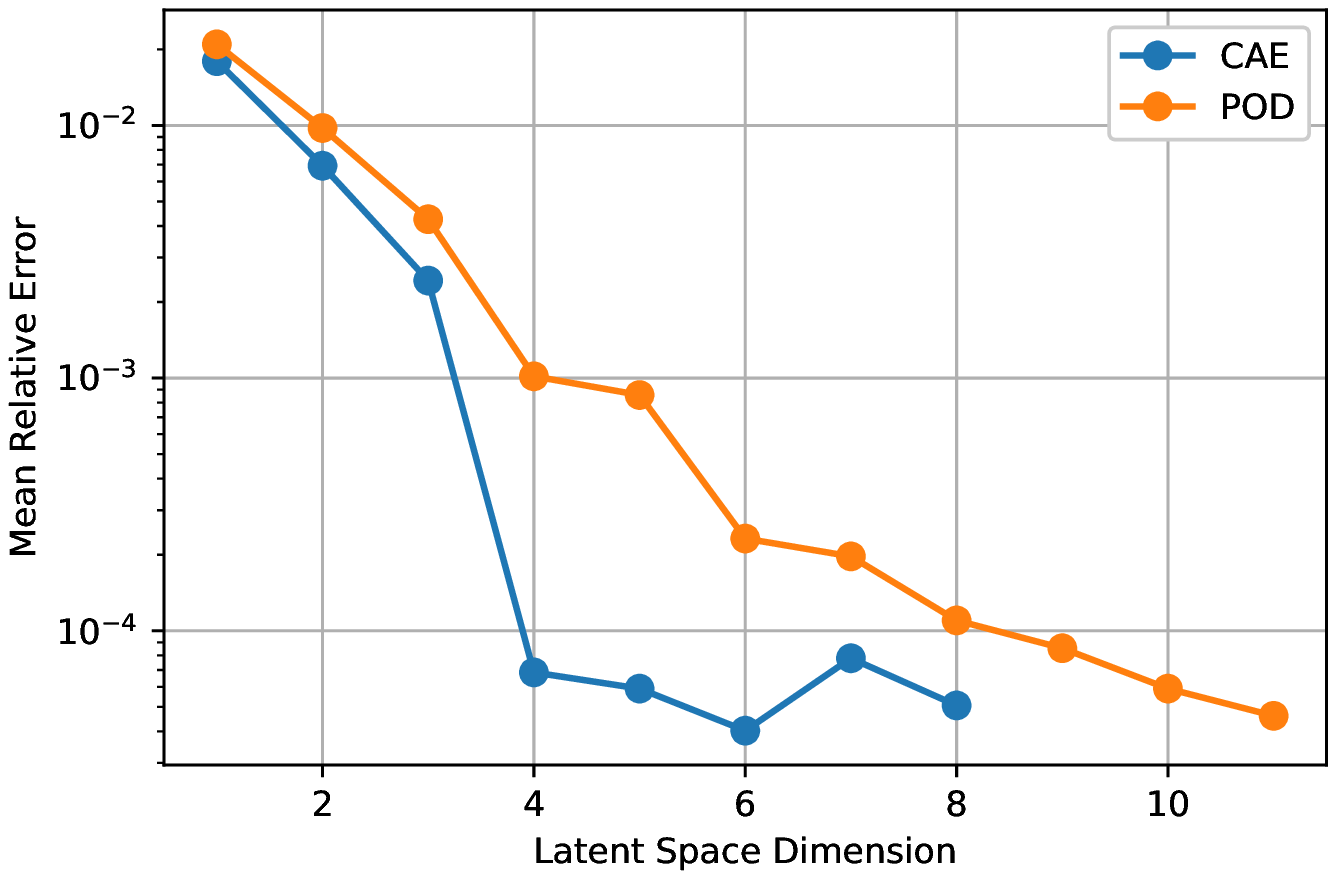}
    \caption{Comparison of convergence of the time averaged MRE of the reconstruction using CAE and POD for the heat equation. }
    \label{fig:POD_AE_comparison}
\end{subfigure}
\caption{}
\end{figure}

In this section we consider results for the linear heat equation parameterized by a space-dependent diffusion rate, $\mu(x)$, on the domain $\Omega=[0,1]^2$: 
\begin{subequations}
\begin{alignat}{2}
    \partial_t u(\mu) - \mu(x)\Delta u(\mu) &= f, &\quad &\text{in} \quad \Omega, \\
    u(\mu) &= u_d, & &\text{on} \quad \Gamma_d, \\
    \textbf{n} \cdot \nabla u(\mu) &= u_n, & &\text{in} \quad \Gamma_n, \\
    u(\mu) &= 0, & &\text{for} \quad t=0.
\end{alignat}
\end{subequations}
We consider a unit square domain divided into four equally-sized square subdomains, $\Omega_1,\Omega_2,\Omega_3,\Omega_4$ with a specific diffusion rate, $\mu_1,\mu_2,\mu_3,\mu_4$ on each subdomain, where $\mu=(\mu_1,\mu_2,\mu_3,\mu_4)\in [0.1,1.5]^4$. Hence, our parameter space is four-dimensional. The boundary conditions are given by $u_d=0$ on $\Gamma_d = \Omega\cap \{y=1\}$, $u_n=0$ on $\Omega\cap\{x=0, x=1\}$, and $u_n=1$ on $\Omega\cap \{y=0\}$  See Figure \ref{fig:heat_setup} for a visualization of the setup.

The high-fidelity snapshots are computed on a $100\times 100$ grid using second-order Lagrange finite elements, which gives third order convergence. The high-fidelity problem has 40401 degrees of freedom. For the implementation we used the FEniCS library in Python \cite{logg2012automated}. The time-stepping is done using the Crank-Nicolson scheme with 100 time steps of size 0.1, meaning the time horizon spans from $t=0$ to $t=10$. 

In Figure \ref{fig:POD_AE_comparison} the convergence of the reconstruction (projection) error for the CAE and the POD is compared. It is clear the CAE achieves higher precision with much lower latent dimension. Indeed, we see that the convergence stagnates around a latent dimension, $N_l=4$, which is actually lower than the intrinsic dimension of the parameterized solution manifold, $\mathcal{M}_{h,\delta t}$. At $N_l=4$ the mean relative error is below $10^{-4}$. To achieve the same accuracy using POD one needs at least a 10-dimensional latent space. The remaining results regarding the heat equation are computed with $N_l=4$.

\begin{figure}
\centering
\begin{subfigure}{.4\textwidth}
  \centering
  \includegraphics[width=.8\linewidth]{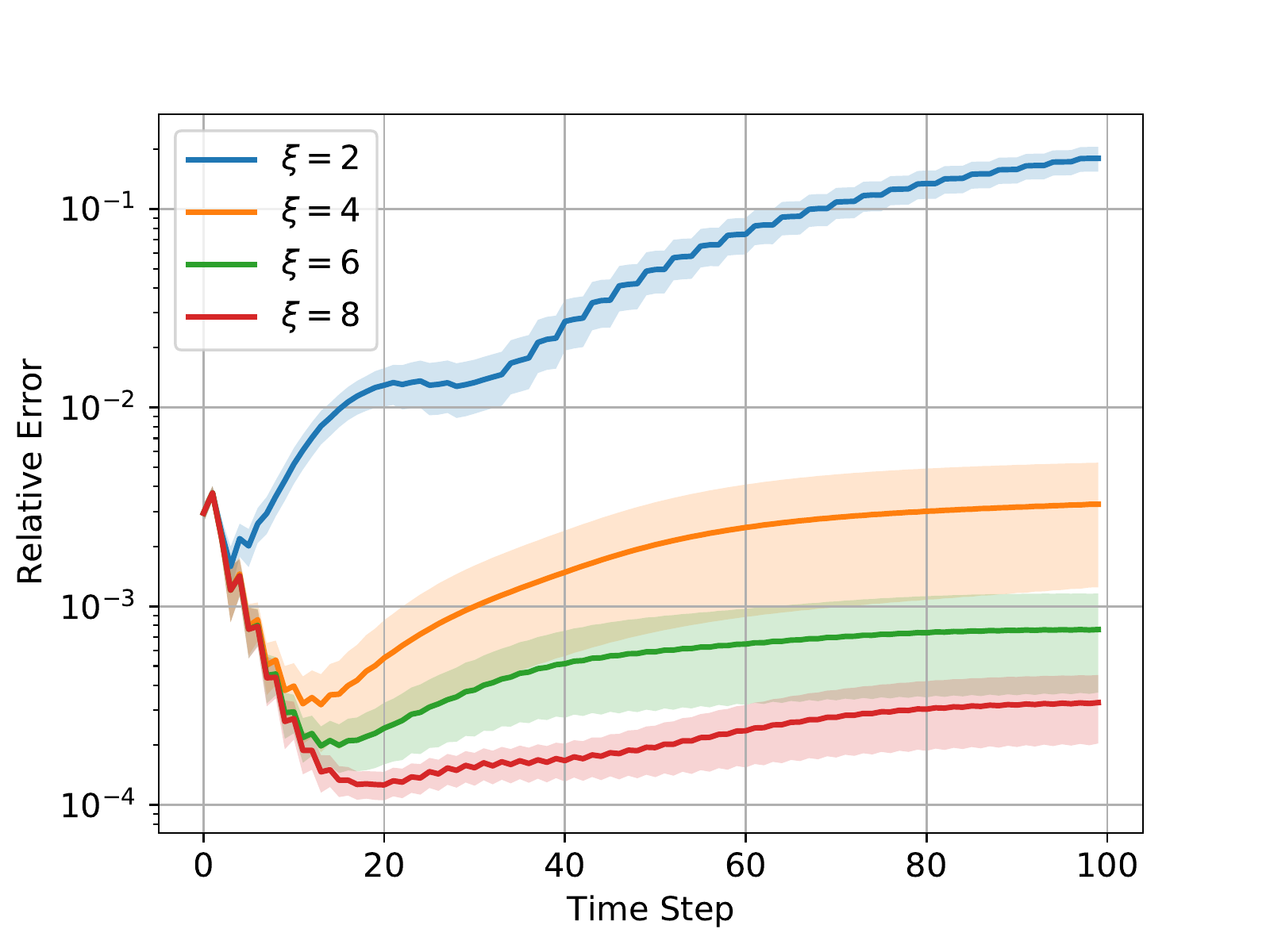}
  \caption{CCNN.}
  \label{fig:heat_CCNN_memory_test}
\end{subfigure}%
\begin{subfigure}{.4\textwidth}
  \centering
  \includegraphics[width=.8\linewidth]{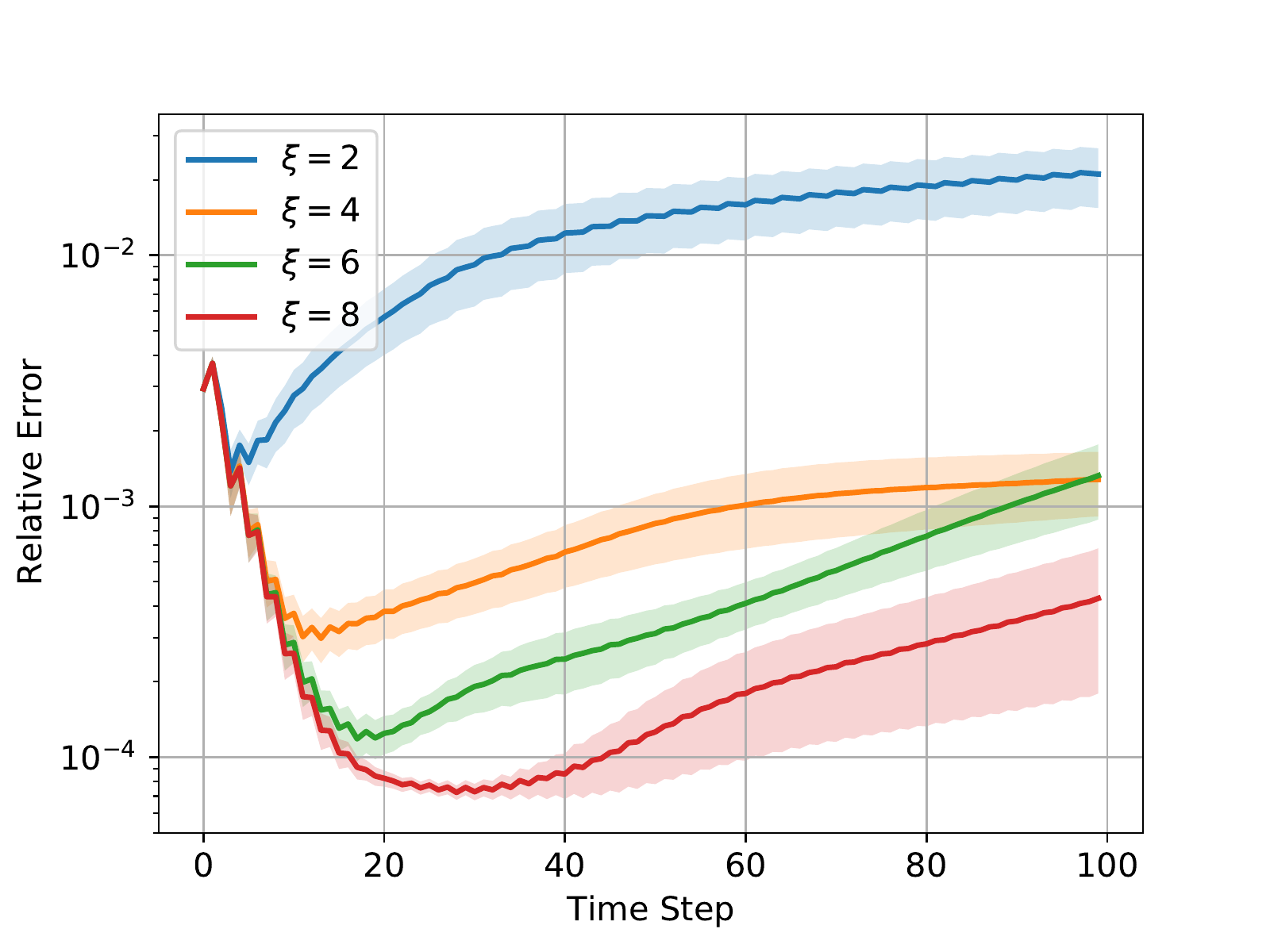}
  \caption{LSTM.}
  \label{fig:heat_LSTM_memory_test}
\end{subfigure}
\caption{Comparison of CCNN and LSTM in relative error for each time step in high-fidelity space for the heat equation for varying memory, $\xi$. The error for each time step is a computed average over 15 test cases with the standard error. }
\label{fig:heat_memory_test}
\end{figure}

In Figure \ref{fig:heat_CCNN_memory_test} and \ref{fig:heat_LSTM_memory_test} we see how the MRE evolves in time for varying memory in the time-stepping NN for the CCNN and LSTM, respectively. The CCNN architecture changes as visualized in Figure \ref{fig:CCNN_types}. As expected, for both the CCNN and the LSTM the error decreases with increasing memory. Furthermore, there is not a big difference between the CCNN and LSTM, but the CCNN seems to perform slightly  better. From hereon, all results regarding the heat equation are computed with a memory of $\xi=8$.

\begin{figure}
\centering
\begin{subfigure}{.4\textwidth}
  \centering
  \includegraphics[width=.8\linewidth]{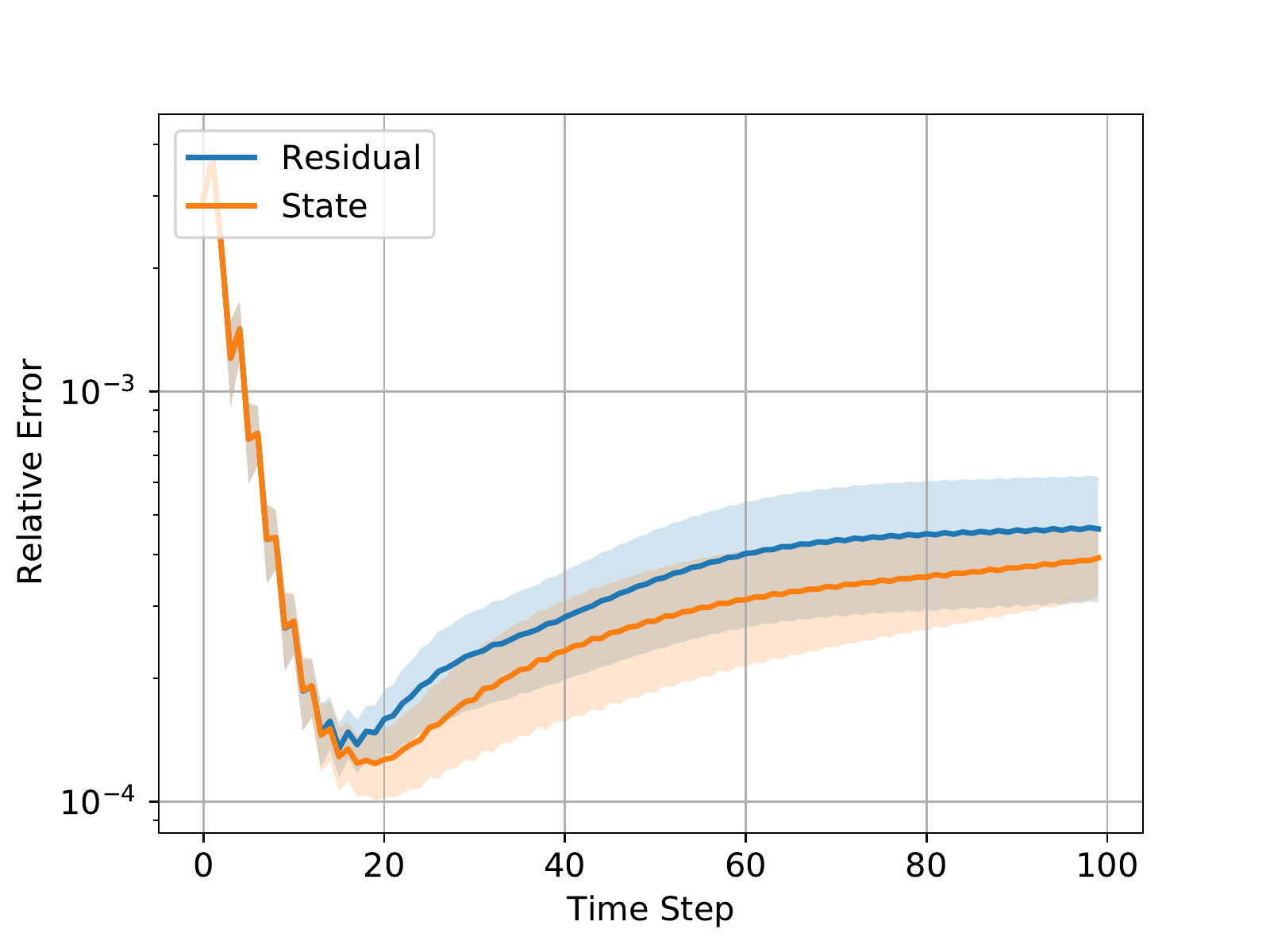}
  \caption{CCNN.}
  \label{fig:heat_CCNN_residual_test}
\end{subfigure}%
\begin{subfigure}{.4\textwidth}
  \centering
  \includegraphics[width=.8\linewidth]{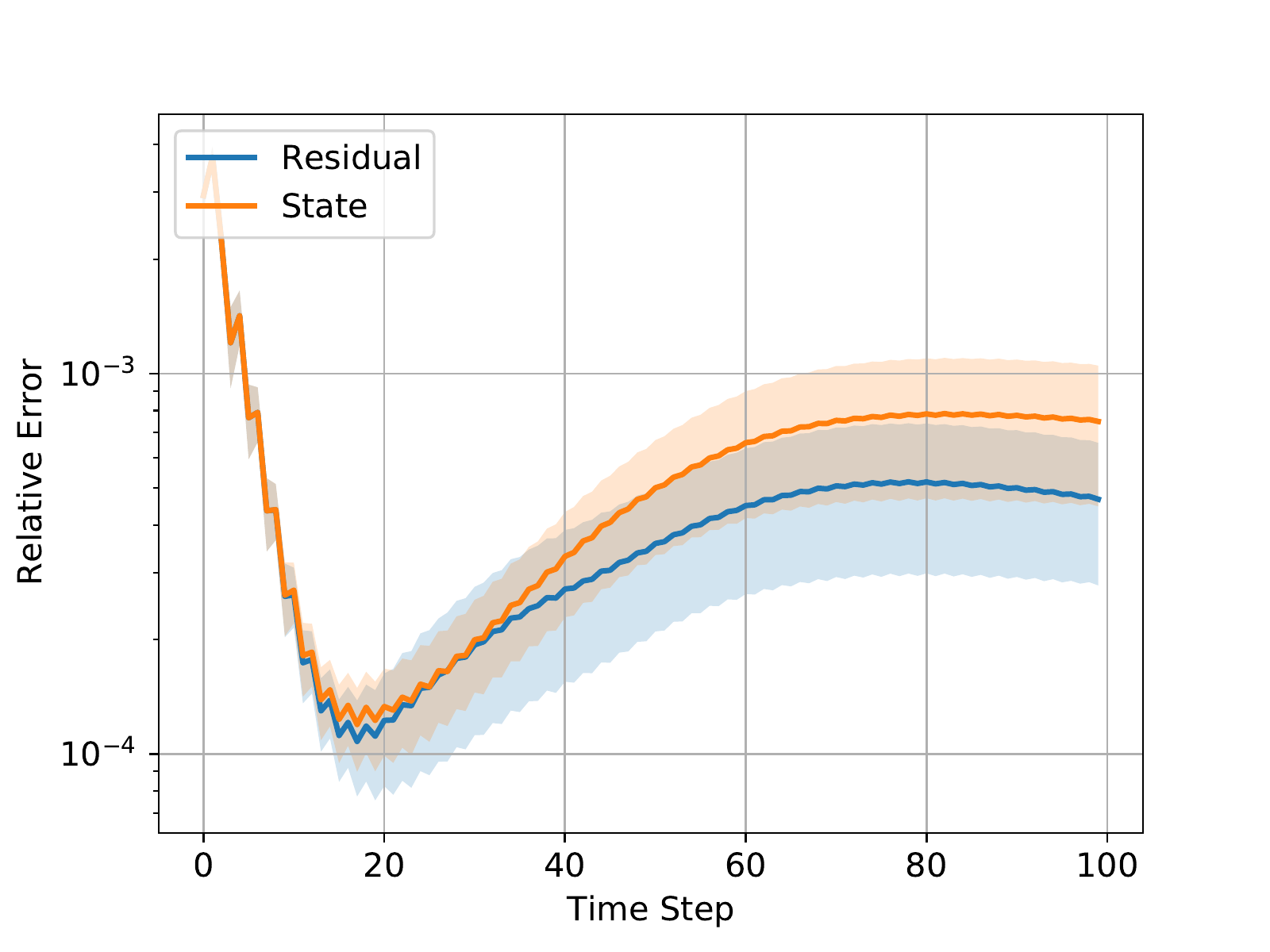}
  \caption{LSTM.}
  \label{fig:heat_LSTM_residual_test}
\end{subfigure}
\caption{Comparison of computing the next step directly and the residual for the CCNN and LSTM. The figures show relative error for each time step in high-fidelity space for the heat equation.}
\label{fig:heat_residual_test}
\end{figure}

One of the proposed improvements of the time-stepping in Section \ref{Approximating Parameterized Time Evolution using Neural Networks} was to approximate the residual instead of the next state directly. In Figure \ref{fig:heat_residual_test} we see that for the CCNN there is close to no difference, with the state approximation performing slightly better. The same results are apparent for the LSTM in Figure \ref{fig:heat_LSTM_residual_test}, but with the residual case being the best choice.  

\begin{figure}
\centering
\begin{subfigure}{.4\textwidth}
  \centering
  \includegraphics[width=.8\linewidth]{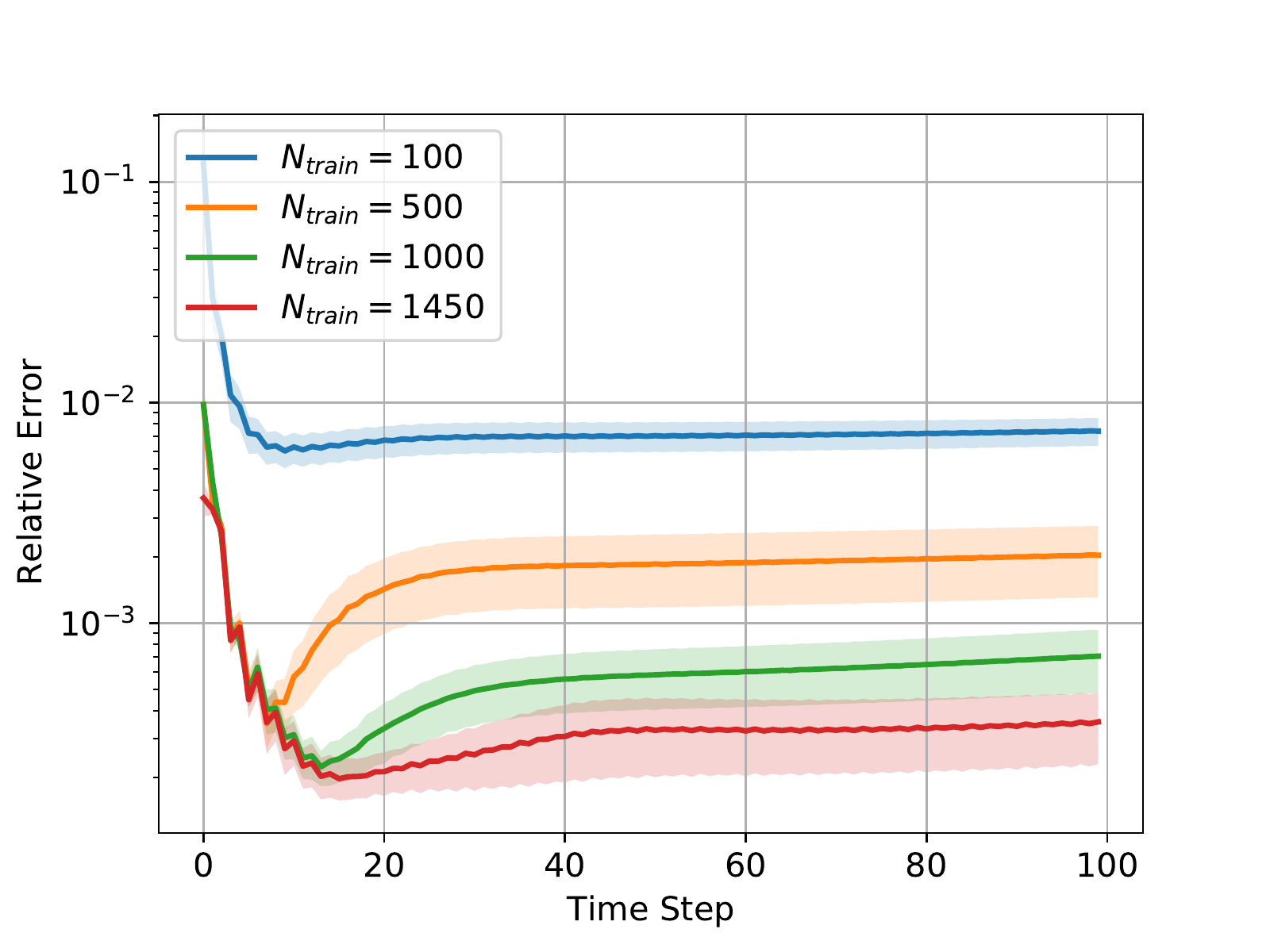}
  \caption{CCNN.}
  \label{fig:heat_CCNN_train_test}
\end{subfigure}%
\begin{subfigure}{.4\textwidth}
  \centering
  \includegraphics[width=.8\linewidth]{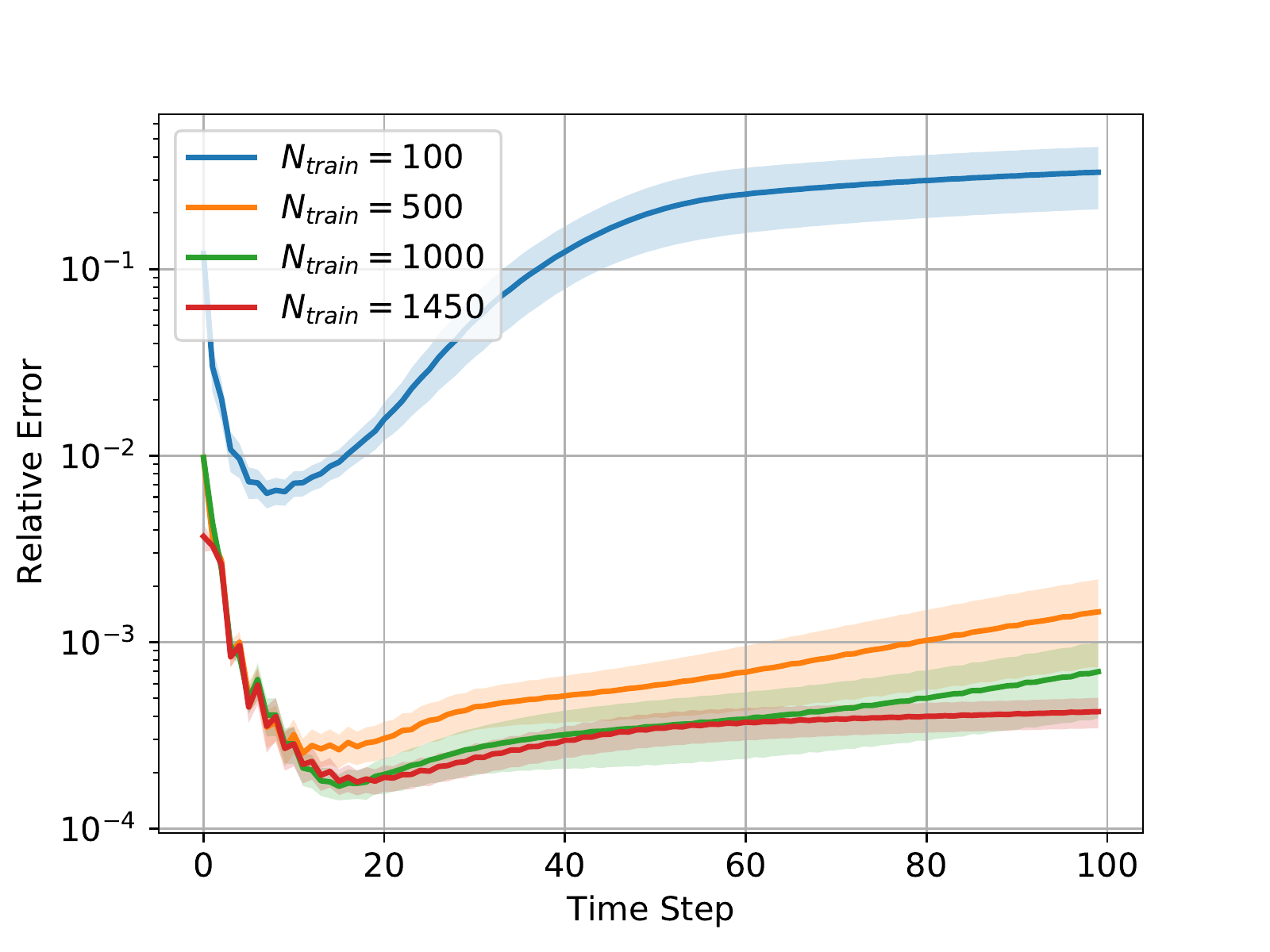}
  \caption{LSTM.}
  \label{fig:heat_LSTM_train_test}
\end{subfigure}
\caption{Comparison of CCNN and LSTM in relative error for each time step in high-fidelity space for the heat equation for number of training samples, $N_{train}$.}
\label{fig:heat_train_test}
\end{figure}

In Figure \ref{fig:heat_CCNN_train_test} and \ref{fig:heat_LSTM_train_test} it is showcased how the number of training trajectories affects the performance for the CCNN and LSTM respectively. As expected, the relative error decreases with increasing number of training trajectories. One further sees that the LSTM requires fewer training trajectories to achieve higher precision in the short term than the CCNN. On the other hand, the MRE seems to stabilize for more time steps in contrast to the LSTM, where the MRE seems to increase. This suggests that the CCNN is the preferable choice for long term predictions with limited training data for this problem. It is worth noting that the relatively large number of training trajectories is due to the fact that we are dealing with a 4-dimensional parameter space, which requires quite a few samples to explorer the full parameter space.

\begin{figure}
\centering
\begin{subfigure}{.33\textwidth}
  \centering
  \includegraphics[width=1.\linewidth]{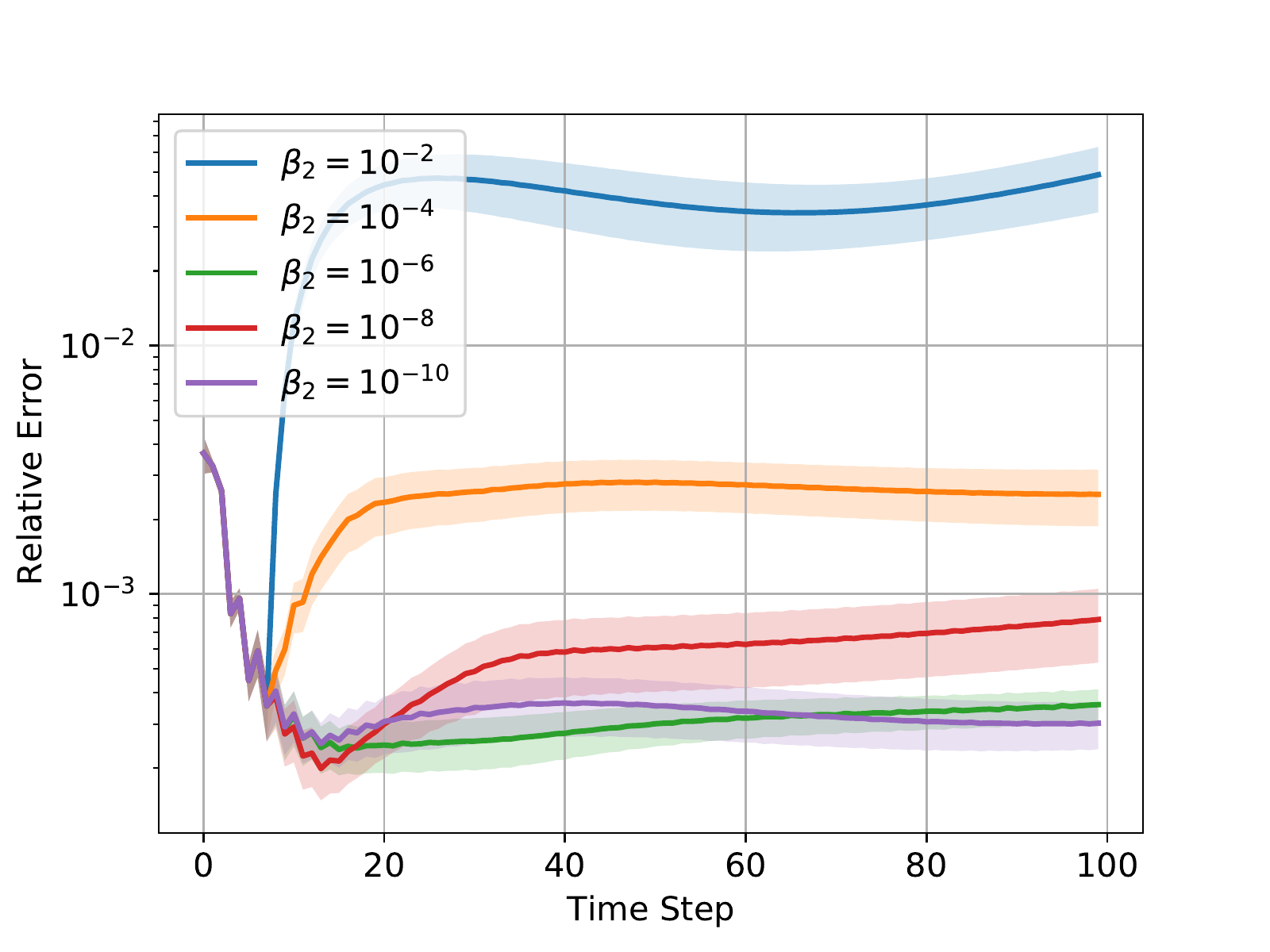}
  \caption{$\beta_1 = 10^{-3}$.}
  \label{fig:heat_CCNN_beta_test_reg3}
\end{subfigure}%
\begin{subfigure}{.33\textwidth}
  \centering
  \includegraphics[width=1.\linewidth]{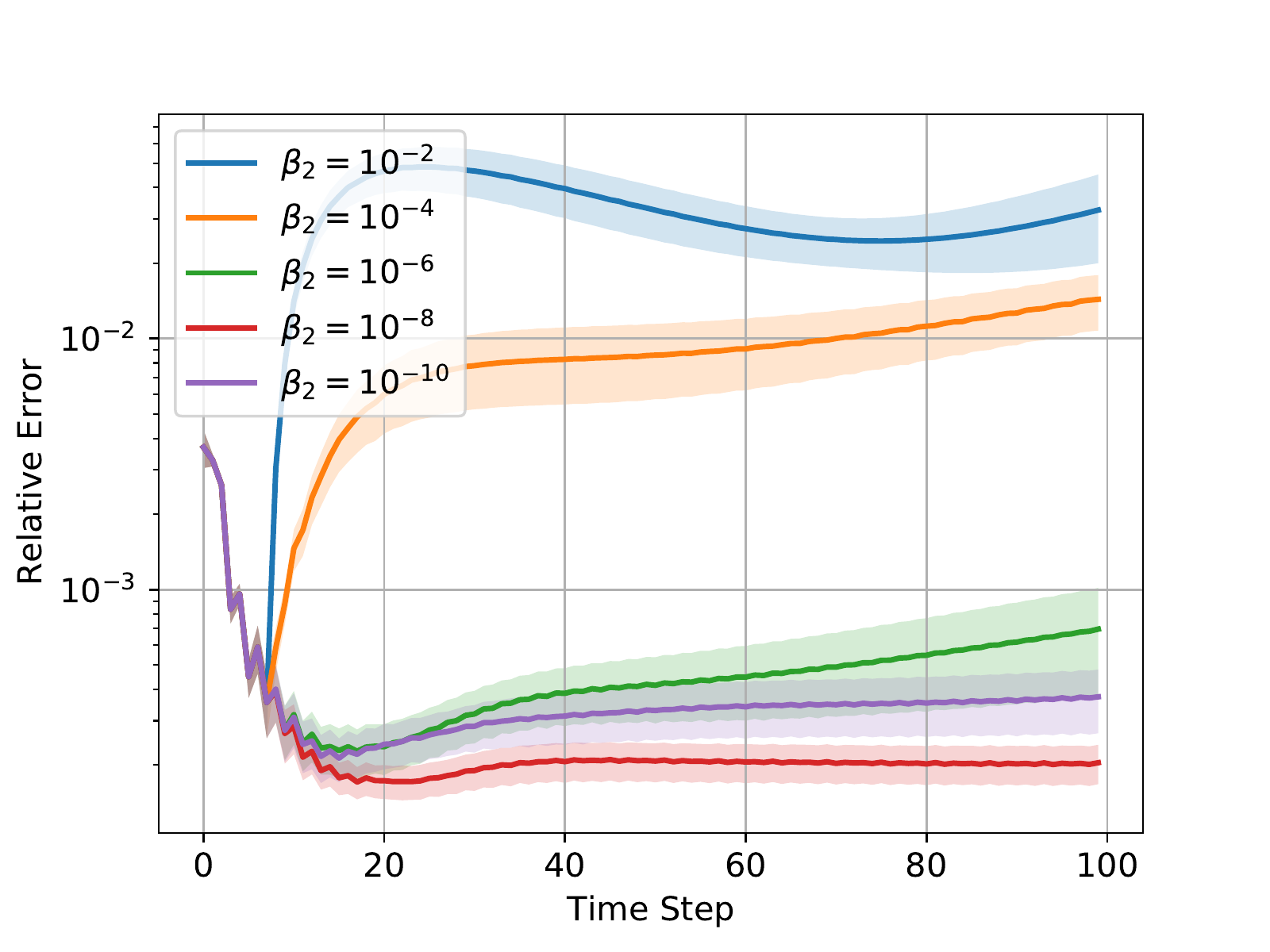}
  \caption{$\beta_1 = 10^{-6}$.}
  \label{fig:heat_CCNN_beta_test_reg6}
\end{subfigure}
\begin{subfigure}{.33\textwidth}
  \centering
  \includegraphics[width=1.\linewidth]{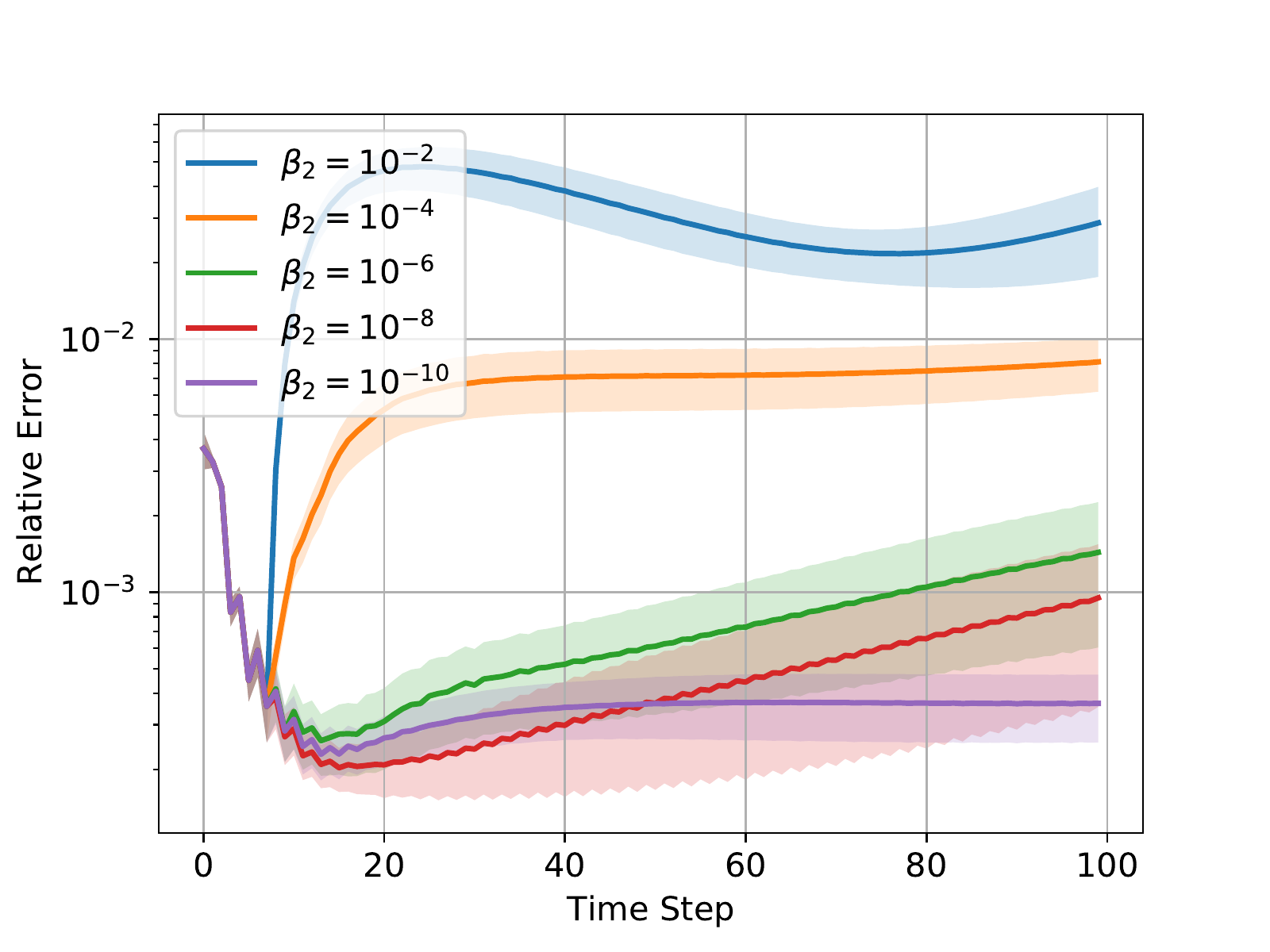}
  \caption{$\beta_1 = 10^{-9}$.}
  \label{fig:heat_CCNN_beta_test_reg9}
\end{subfigure}
\caption{Impact of the two regularization terms, weight decay, $\beta_1$, and Jacobian, $\beta_2$ for the heat equation. The average relative error in high-fidelity space over 15 test trajectories for each time step is shown. Each figure shows the error for a constant $\beta_1$ and varying $\beta_2$. }
\label{fig:heat_beta_test}
\end{figure}

In Section \ref{Approximating Parameterized Time Evolution using Neural Networks} we discussed two regularization terms, weight decay, $\beta_1$, and Jacobian, $\beta_2$, with the aim of promoting generalization and stability. The effect of those terms is shown in Figure \ref{fig:heat_beta_test} for the CCNN. It is, somewhat surprisingly, apparent that the size of the weight decay regularization term has very little effect on the performance, while it is clear that the Jacobian regularization has a large impact. In general, the performance improves with smaller amount of regularization suggesting that stability is not an issue. This could be because the dynamics in this diffusion type problem rather quickly go to a steady state.

\begin{figure}
\centering
\begin{subfigure}{.45\textwidth}
  \centering
  \includegraphics[width=.8\linewidth]{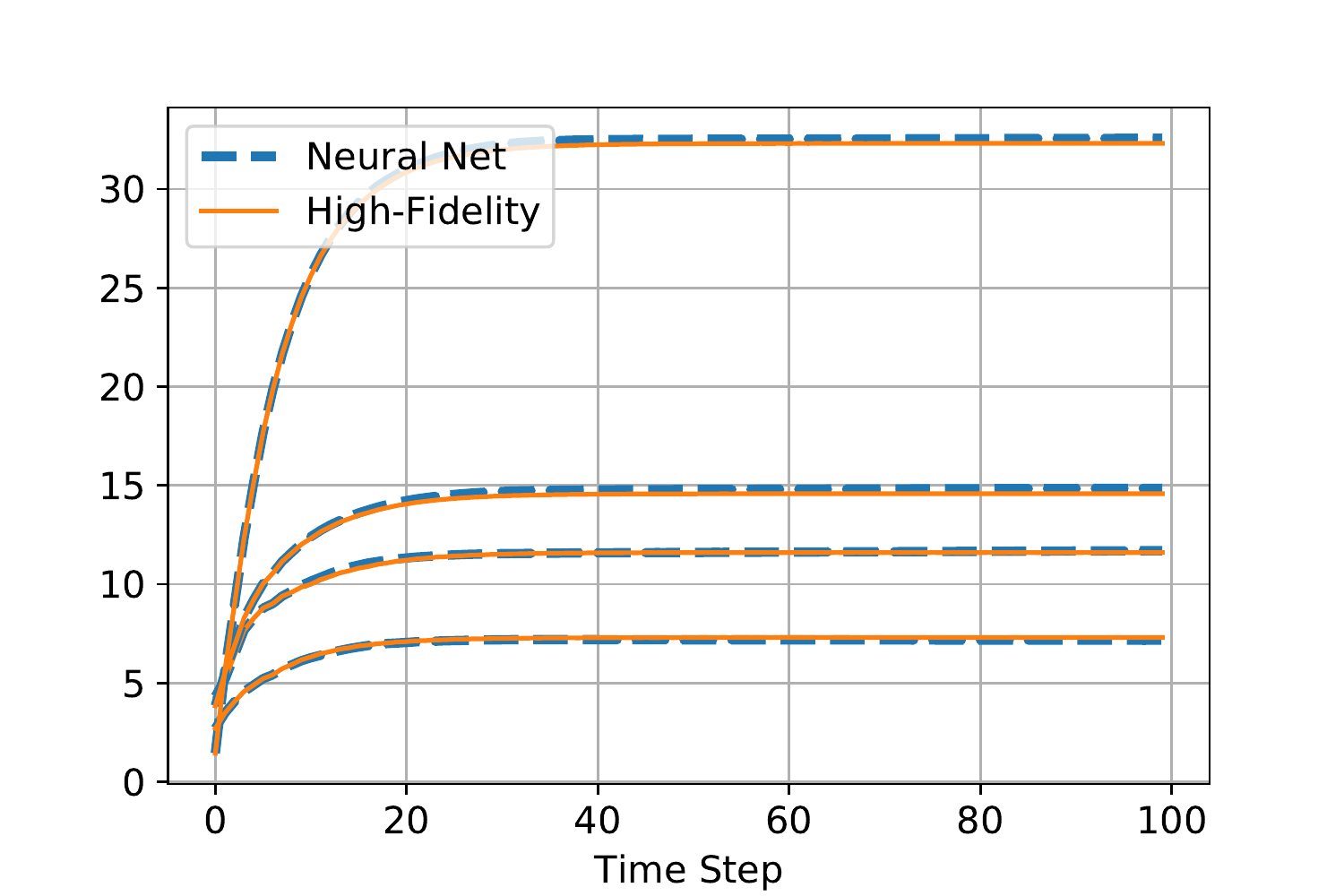}
  \caption{Latent space trajectories for the heat equation.}
  \label{fig:reduced_heat_sols}
\end{subfigure}%
\begin{subfigure}{.45\textwidth}
  \centering
  \includegraphics[width=.8\linewidth]{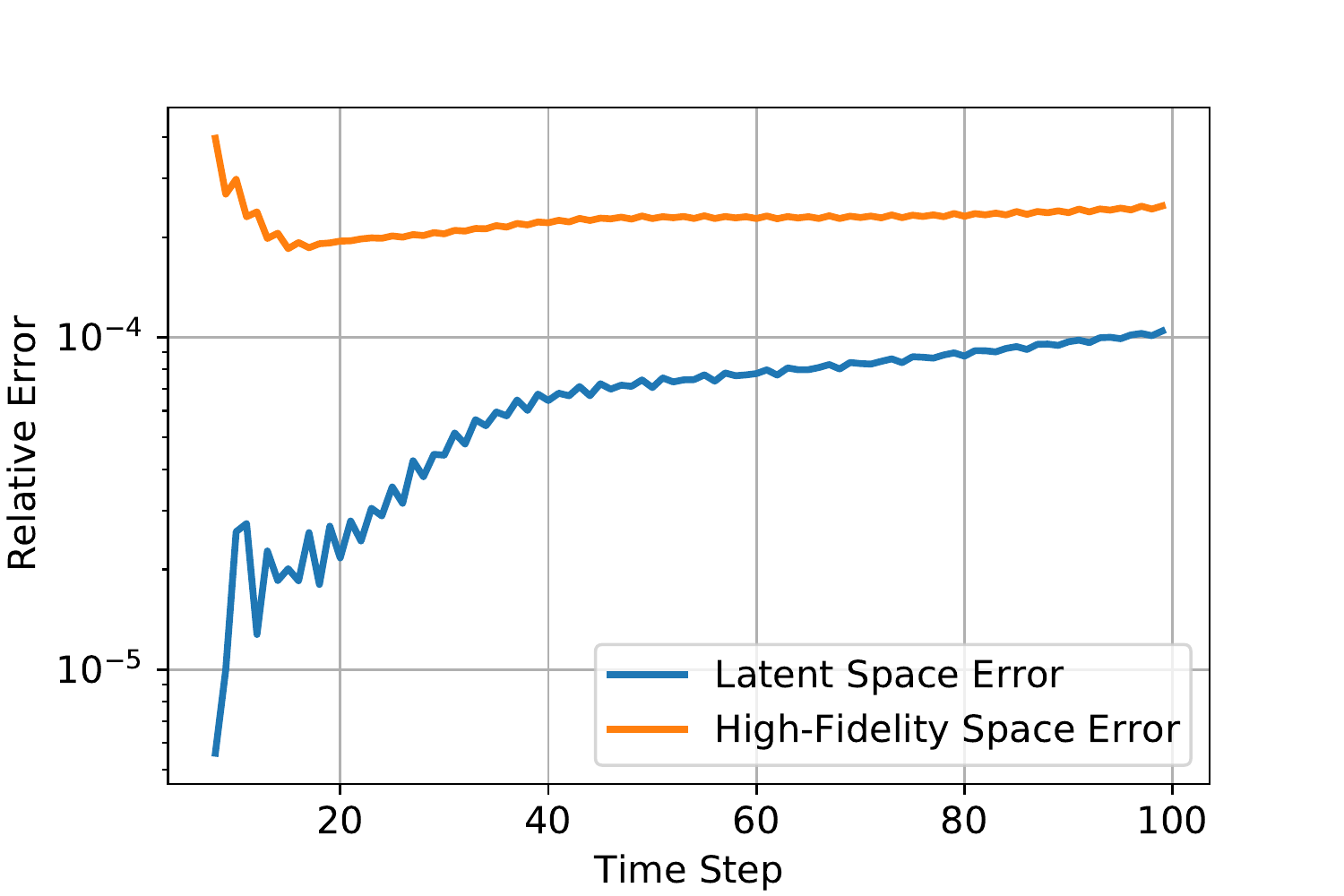}
  \caption{Mean relative error in the latent space and high-fidelity space for the heat equation.}
  \label{fig:heat_time_error_CCNN}
\end{subfigure}
\caption{$\mu =(0.1902, 0.4583 , 1.2648, 0.7116)$.}
\label{fig:heat_latent_sol_and_error}
\end{figure}

In Figure \ref{fig:reduced_heat_sols} we see the high-fidelity trajectories reduced to the latent space as well as the NN time-stepping approximation for a specific test case. It is clear that in the latent space the dynamics reach a steady state from around time step 20. Furthermore, the approximated trajectories are very close to the reduced high-fidelity trajectories. In figure \ref{fig:heat_time_error_CCNN} a comparison of the error in the latent space,  $\text{MRE}(u_l(\mu),\tilde{u}_l(\mu))$, and in the high-fidelity space, $\text{MRE}(u_h(\mu),\tilde{u}_h(\mu))$, is made. Throughout all time steps the error in the latent space is lower than in the high-fidelity space. This increase in error comes from the decoding step of the CAE. 

\begin{figure}
\centering
\begin{subfigure}{.33\textwidth}
  \centering
  \includegraphics[width=1.\linewidth]{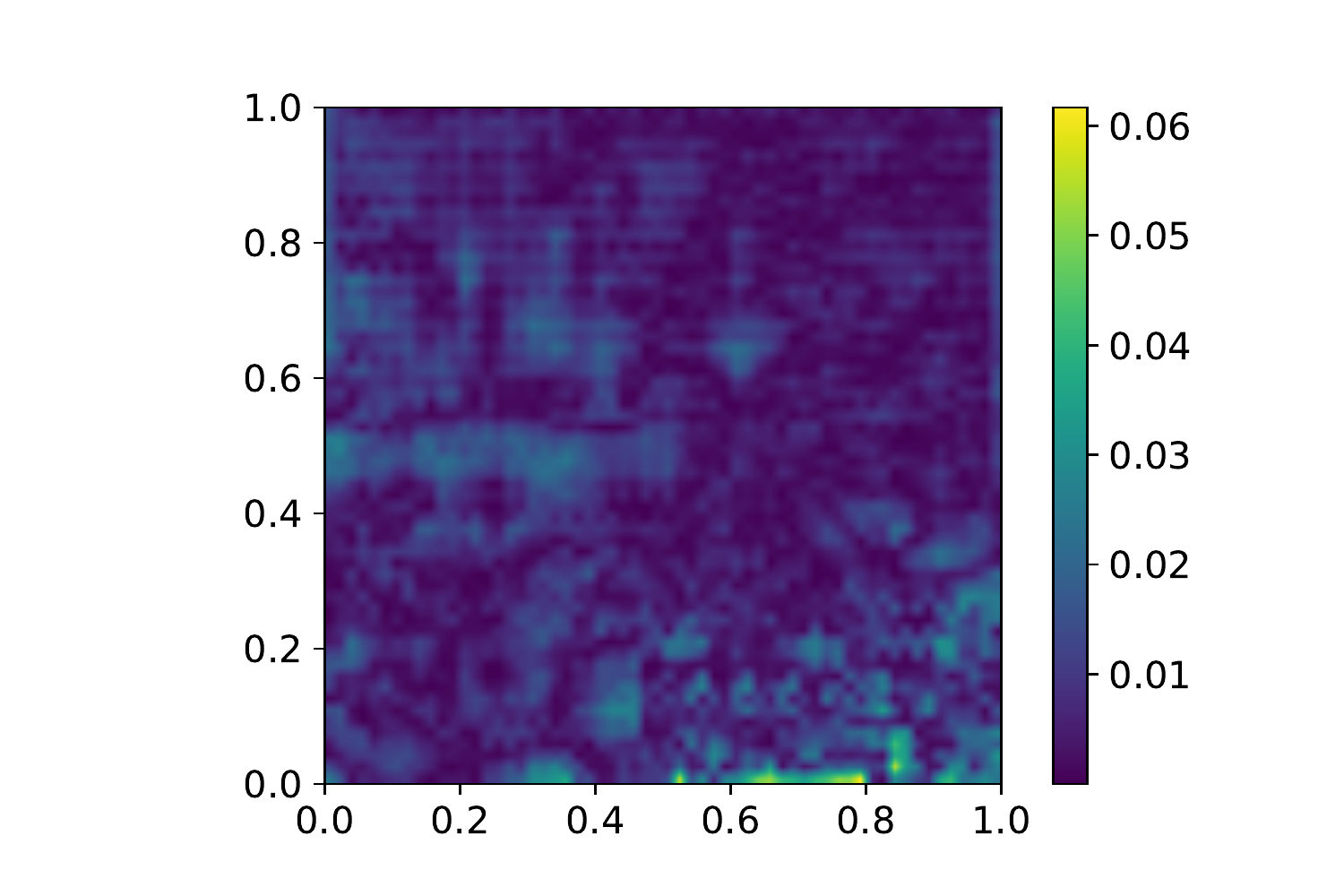}
  \caption{$t=0$}
  \label{fig:heat_error_t0}
\end{subfigure}
\begin{subfigure}{.33\textwidth}
  \centering
  \includegraphics[width=1.\linewidth]{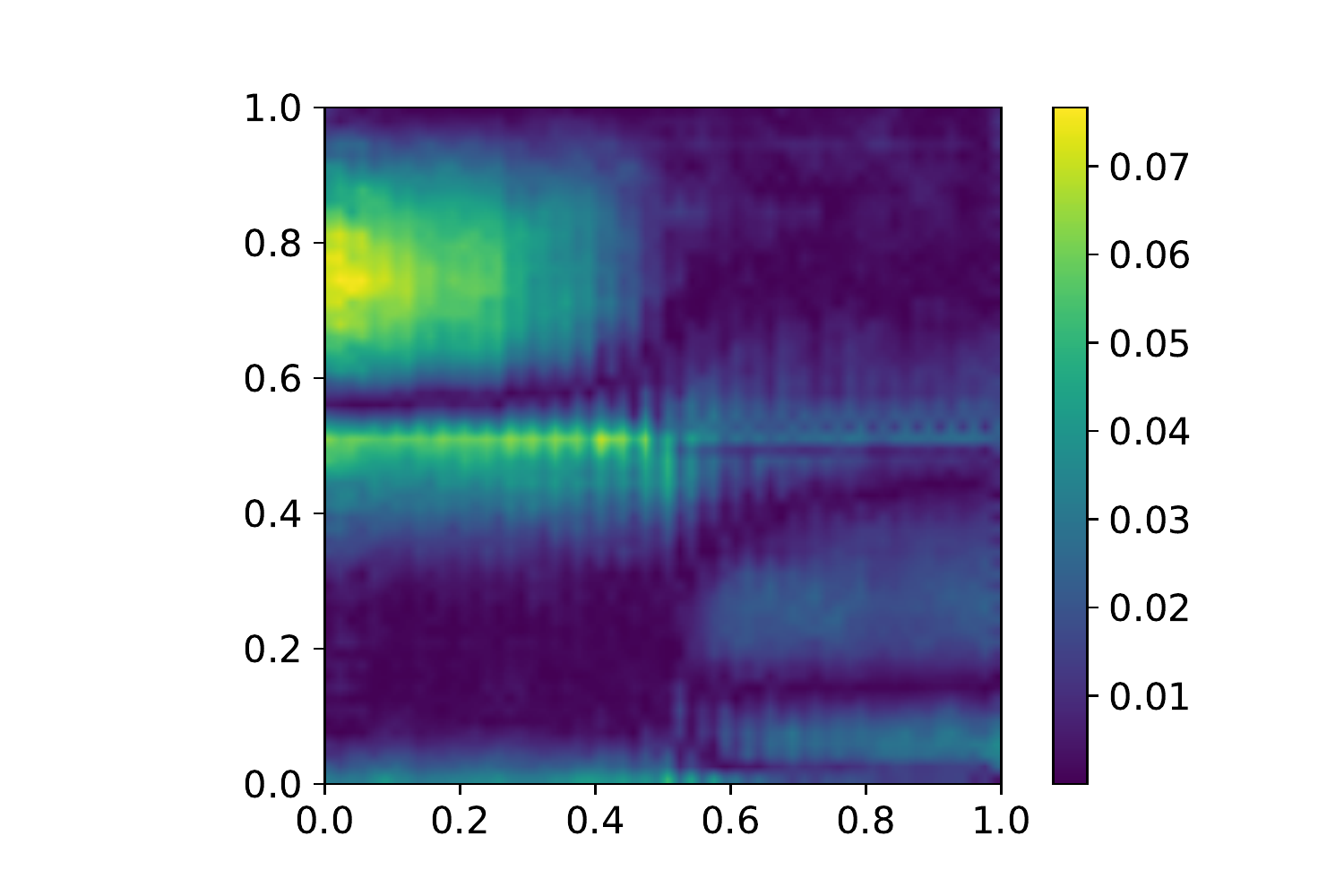}
  \caption{$t=1$}
  \label{fig:heat_error_t1}
\end{subfigure}%
\begin{subfigure}{.33\textwidth}
  \centering
  \includegraphics[width=1.\linewidth]{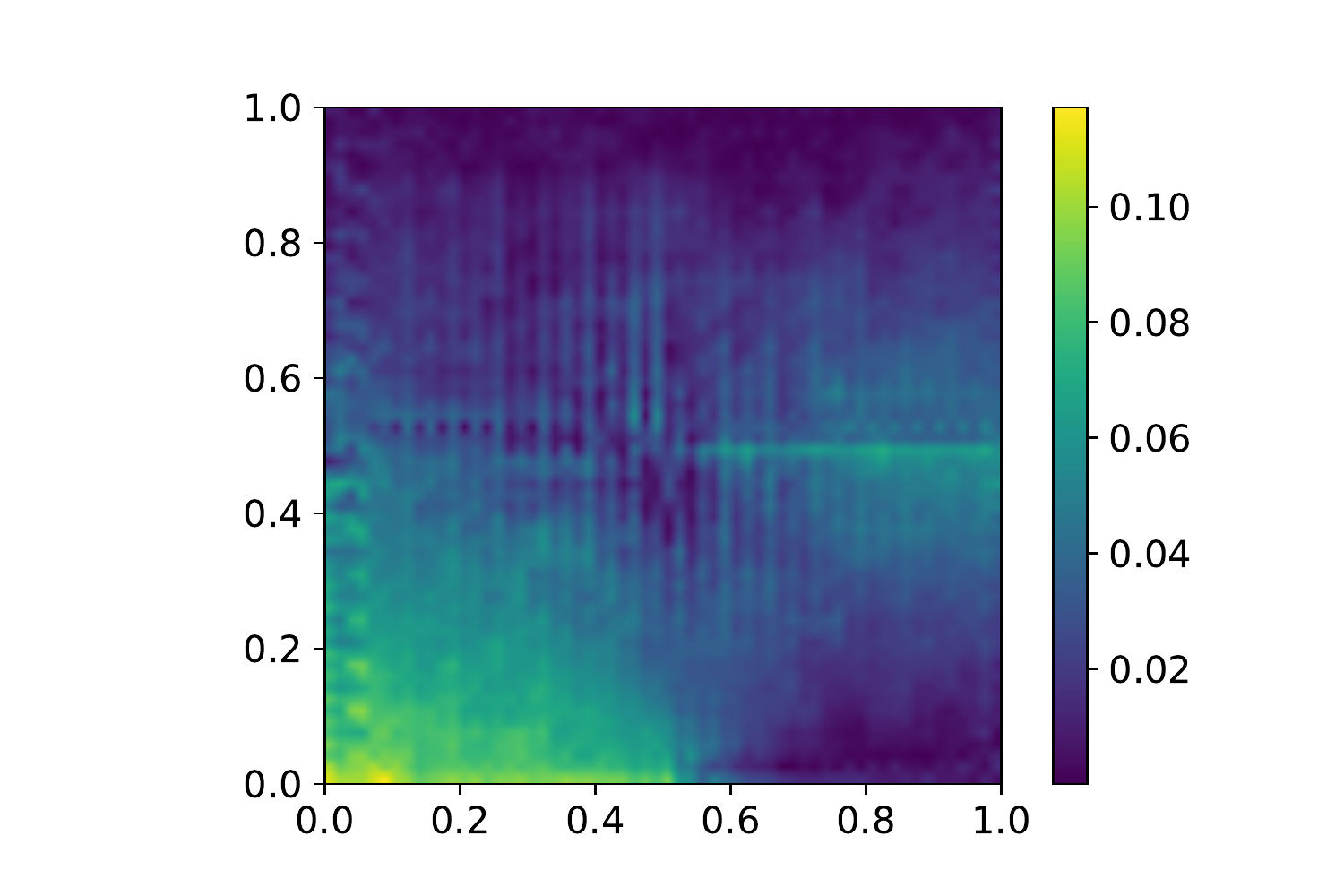}
  \caption{$t=10$}
  \label{fig:heat_error_t10}
\end{subfigure}
\caption{Pointwise absolute error between the high-fidelity solution and the neural network prediction (CCNN) for the heat equation with diffusion rates, $\mu=(1.0904, 0.3550, 0.1267, 1.0159 )$ .}
\label{fig:heat_error}
\end{figure}

For time steps $t=0$, $t=1$, and $t=10$, the pointwise absolute error is shown for a specific test case in Figure \ref{fig:heat_error}.

\subsection{Linear Advection Equation }
We consider a linear advection equation on the domain $\Omega=[0,1]^2$:
\begin{subequations}
\begin{alignat}{2}
   \partial_t u(\mu) + b\cdot \nabla u(\mu) &= 0, & \quad &\text{in} \quad \Omega, \\
   u(\mu) &= 0 & & \text{on} \quad \Gamma,
\end{alignat}
where $\Gamma=\partial\Omega$,  
\begin{align}
    b = \mu_1 \begin{pmatrix}
        - y - \frac{1}{2} \\ x-\frac{1}{2}
    \end{pmatrix},
\end{align}
\end{subequations}
with initial condition
\begin{align}
u_0(\mu) = \exp\left(\frac{1}{2}\left[\frac{(x-x_{0})^2}{0.005} + \frac{(y-y_{0})^2}{0.005}\right]\right),
\end{align}
where
\begin{align}
    \begin{pmatrix} x_0 \\ y_0 \end{pmatrix} = 
    \frac{1}{4}\begin{pmatrix} \cos(\mu_2)  \\ \sin(\mu_2)  \end{pmatrix}+ \frac{1}{2}.
\end{align}
This problem models a Gaussian curve being advected with velocity $\mu_1$ in a circle with origin at $[\frac{1}{2},\frac{1}{2}]$ and radius $\frac{1}{4}$, starting at the position given by the angle $\mu_2$. This problem is parameterized by two parameters, $\mu=(\mu_1,\mu_2)\in [0.5,1.5]\times [0,2\pi]$. The first parameter, the velocity, is directly affecting the dynamics, while the other, $\mu_2$, is only dictating the initial placement of the Gaussian curve. Hence, we are dealing with a 2-dimensional parameter space, while the dynamics are only parameterized by a single parameter. 

The high-fidelity snapshots are computed on a $60\times 60$ grid using the discontinuous Galerkin method with linear Lagrange elements, resulting in a second-order convergence scheme that suits advection dominated problems well. The high-fidelity model consists of 21600 degrees of freedom. For the implementation we used the FEniCS library in Python \cite{logg2012automated}. The time-stepping is done using the Crank-Nicolson scheme with time steps of size 0.0075 for 2000 steps. resulting in a time interval, $t\in[0,15]$. The training of the neural networks is done using every 4th time step, $s=4$, meaning the model is trained to take steps of size 0.03.  

The same tests as for the heat equation have been performed on the advection equation. The NN configuration that performs the best uses a memory of $\xi=6$, $\beta_1=10^{-9}$, $beta_2=10^{-6}$, and computes the residual rather than the state directly. Furthermore, the training has been performed with 120 training trajectories.

\begin{figure}
\centering
\begin{subfigure}{.33\textwidth}
    \centering
    \includegraphics[width=1.\textwidth]{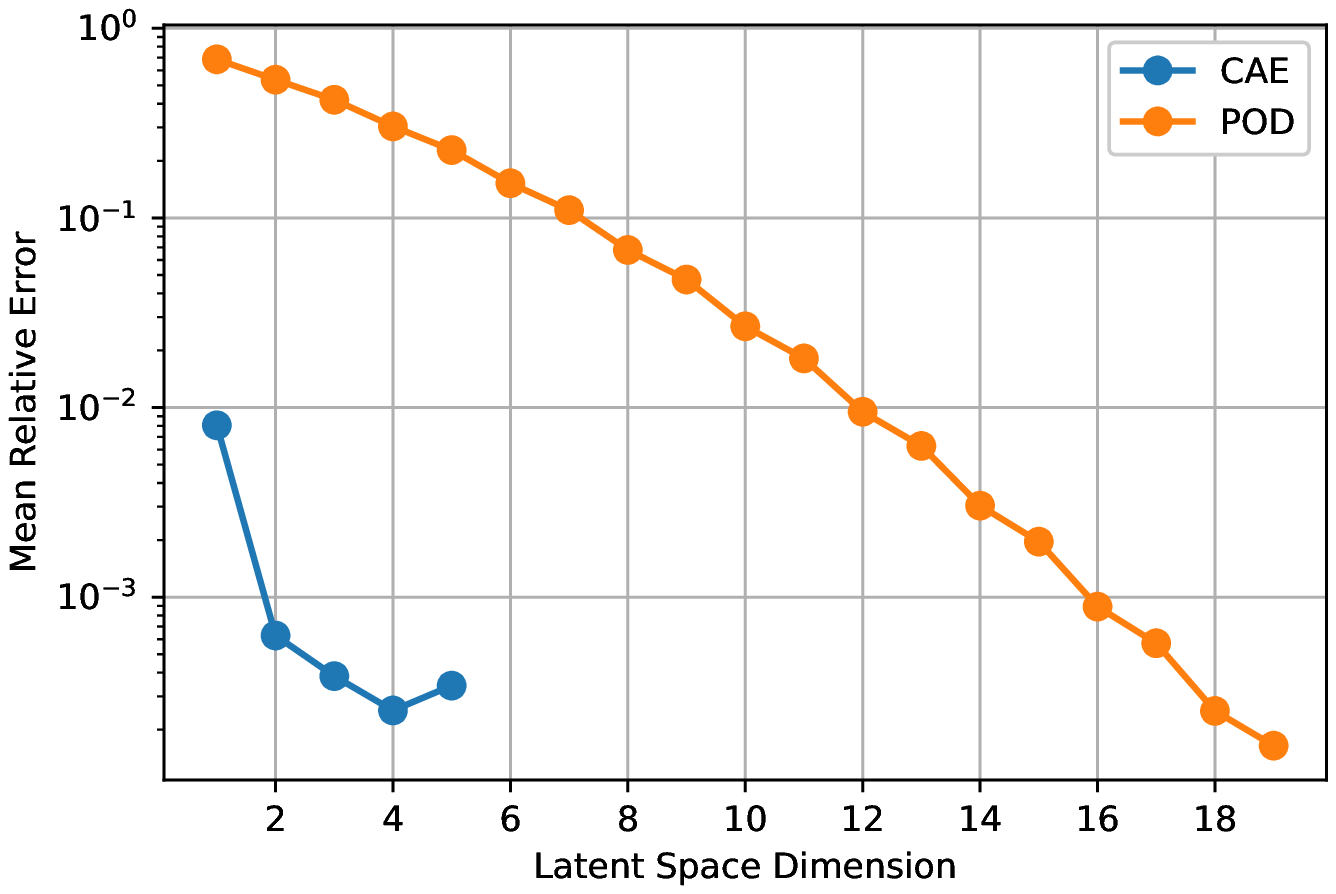}
    \caption{}
    \label{fig:advection_POD_vs_AE}
\end{subfigure}
\begin{subfigure}{.33\textwidth}
  \centering
  \includegraphics[width=1.\linewidth]{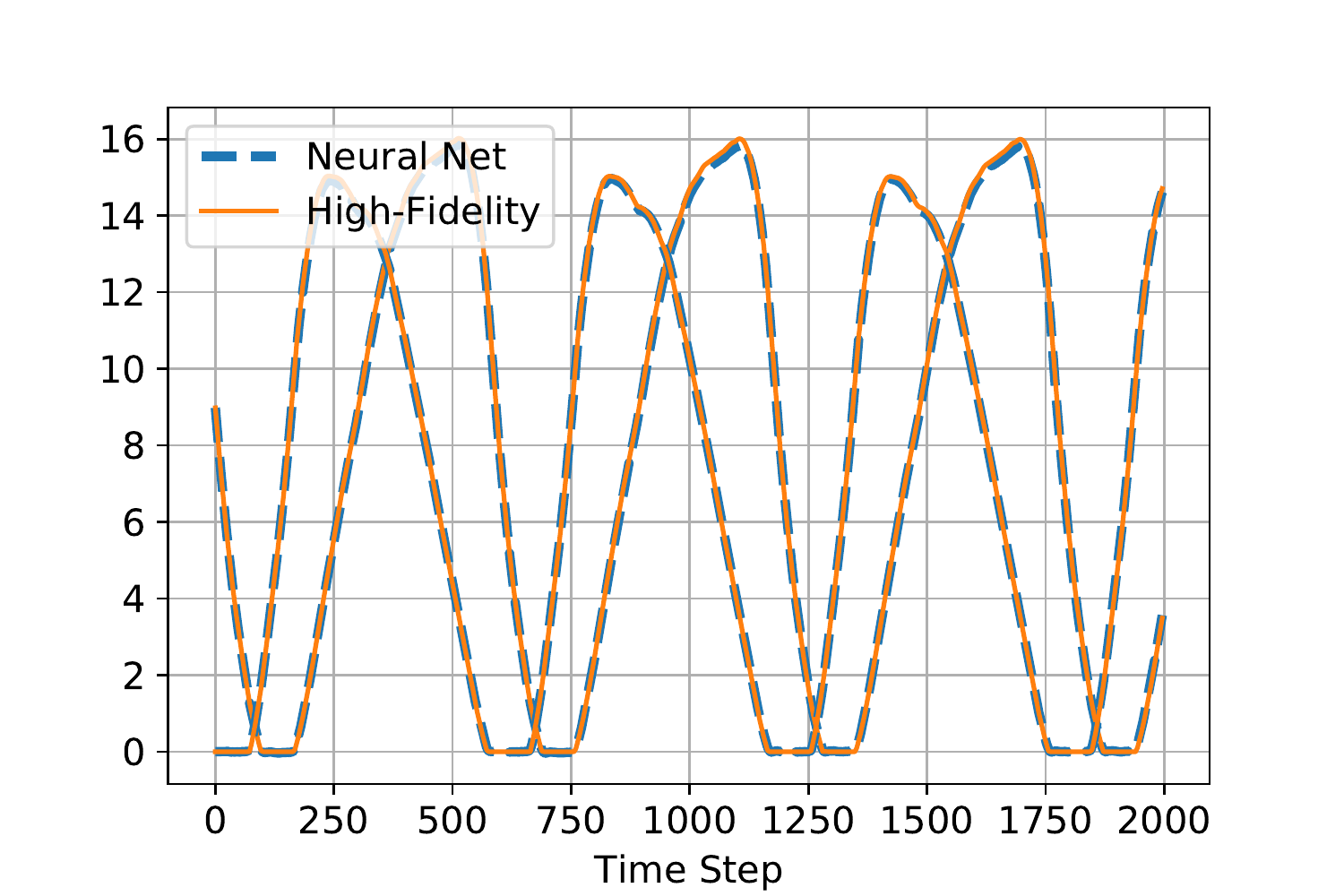}
  \caption{}
  \label{fig:reduced_adv_sols}
\end{subfigure}%
\begin{subfigure}{.33\textwidth}
  \centering
  \includegraphics[width=1.\linewidth]{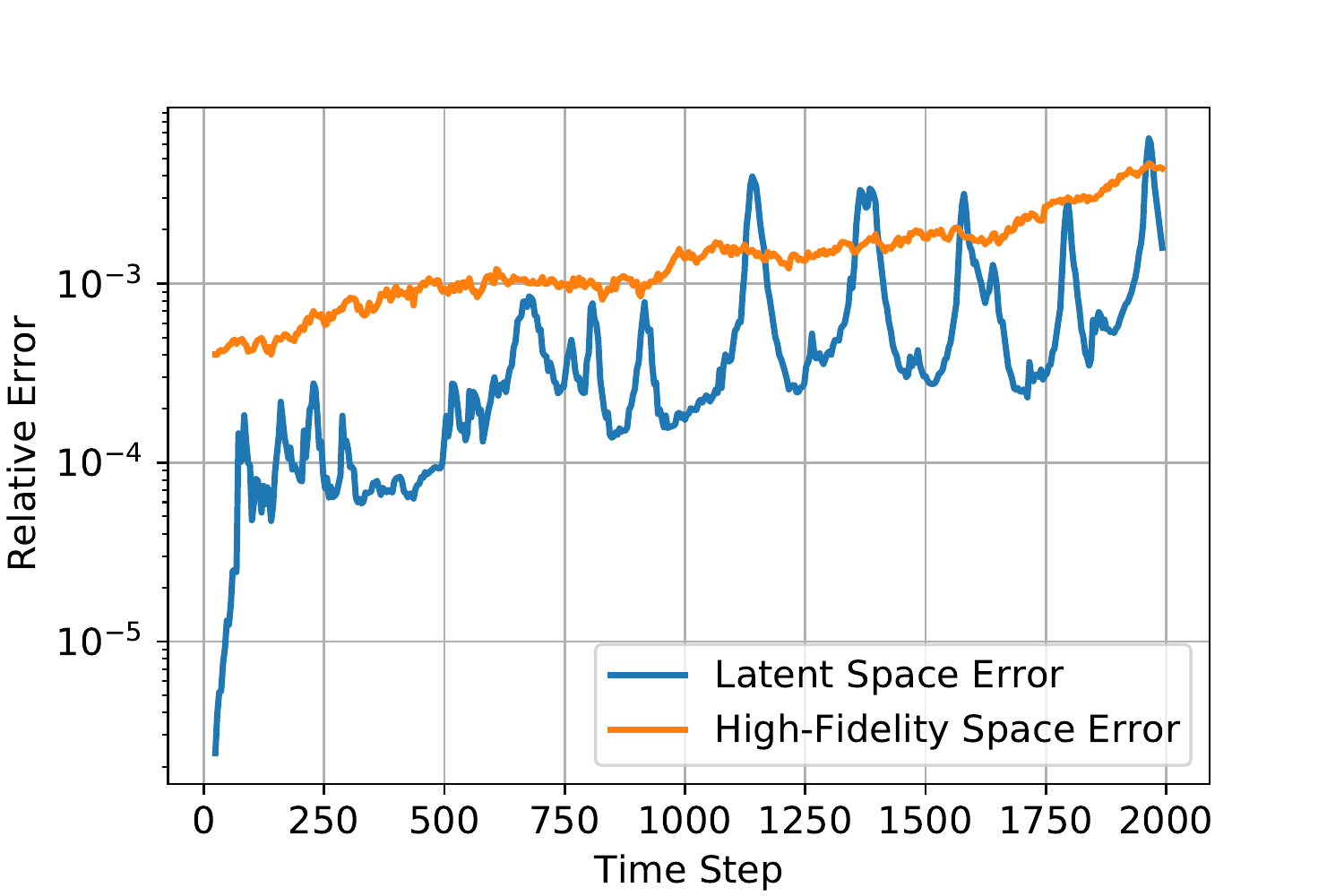}
  \caption{}
  \label{fig:adv_time_error_CCNN}
\end{subfigure}
\caption{(a) Comparison of convergence of the time averaged MRE of the reconstruction using CAE and POD for the heat equation. (b) Latent space trajectories with advection velocity $\mu_1=1.4161$, and initial angle, $\mu_2=2.8744$ and (c) average test errors computed for the linear advection equation. The predictions are computed using the CCNN dynamic network. }
\label{fig:adv_reduced_sols}
\end{figure}

In Figure \ref{fig:advection_POD_vs_AE}, we see a significant improvement by using the CAE compared to the POD approach. Using a latent dimension of 2, which is also the intrinsic dimension of the solution manifold, the CAE reconstructs the high-fidelity solution with an MRE between $10^{-3}$ and $10^{-4}$. To achieve the same accuracy using the POD method, one needs a latent dimension of at least 17. This supports the previous claim that POD does, in general, not perform well on advection dominated problems. 

Figure \ref{fig:reduced_adv_sols} shows the encoded high-fidelity trajectories as well as the latent space trajectories computed by using the NN. As expected, we see periodic behavior and close to no discrepancy between the encoded high-fidelity trajectory and the NN computed latent space trajectory. 

Comparing the errors in the latent space with the errors in the high-fidelity space in Figure \ref{fig:adv_time_error_CCNN}, we observe that the latent space errors are in general an order of magnitude smaller. However, the errors in the latent space show much more spurious oscillations, suggesting that the CAE is not sensitive to small perturbations in the latent space. 

\begin{figure}
\centering
\begin{subfigure}{.33\textwidth}
  \centering
  \includegraphics[width=1.\linewidth]{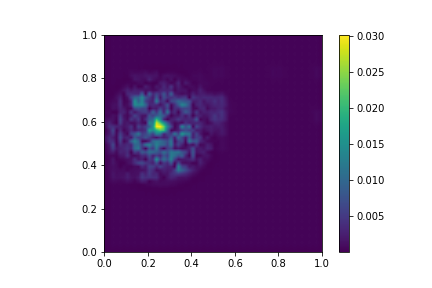}
  \caption{$t=0$}
  \label{fig:linear_adv_error_t0}
\end{subfigure}
\begin{subfigure}{.33\textwidth}
  \centering
  \includegraphics[width=1.\linewidth]{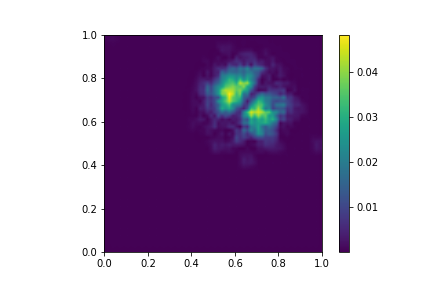}
  \caption{$t=7.5$}
  \label{fig:linear_adv_error_t75}
\end{subfigure}%
\begin{subfigure}{.33\textwidth}
  \centering
  \includegraphics[width=1.\linewidth]{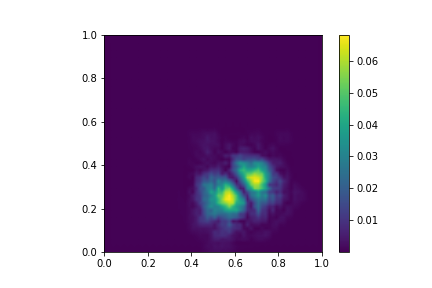}
  \caption{$t=15$}
  \label{fig:linear_adv_error_t15}
\end{subfigure}
\caption{Pointwise absolute error between the high-fidelity solution and the neural network prediction (CCNN) for the linear advection equation with velocity, $\mu_1=1.41605949$, and initial angle, $\mu_2=2.87439799$.}
\label{fig:linear_advection_sols}
\end{figure}

By looking at the pointwise error between the high-fidelity and the NN solutions in Figure \ref{fig:linear_advection_sols} it is clear that the NN approximation introduces a small phase error. This error is, however, not large and could possibly be corrected in a post-processing step. 

\subsection{2D Nonlinear Equation - Lid Driven Cavity}
In this section we showcase results for the lid driven cavity incompressible Navier-Stokes problem, parameterized by the Reynolds number:

\begin{subequations}
\begin{alignat}{3}
    \partial_t u(Re) + (u(Re)\cdot \nabla)u(Re) - \nabla p(Re) &= \frac{1}{Re} \Delta u(Re),  &\text{in} \quad \Omega, \\ 
    \nabla \cdot u(Re) &= 0,  &\text{in} \quad \Omega, \\
    u(Re) &= 0, &\text{on} \quad \Gamma, \\
    (u_x(Re),u_y(Re)) &= (1,0) &\text{for} \quad y=0,
\end{alignat}
\end{subequations}
where $\Omega = (0,1)\times (0,1)$ and $\Gamma$ is the boundary except for the one with $y=0$, and the initial condition is $u(Re)=0$ and $p(Re)$. The goal is to compute the velocity field, $u(Re) = (u_x(Re),u_y(Re))$ and the pressure field, $p(Re)$ in the full time domain and for all Reynolds numbers $Re$ in an interval. Thus, $\mu=Re$. We consider $Re\in[100,300]$. The high-fidelity snapshots are computed using a second-order finite volume scheme on a $100\times 100$ staggered non-uniform grid, ensuring higer resolution near the boundaries. This gives $3\cdot100\cdot100=30000$ degrees of freedom - the $x$-component and the $y$-component of the velocity field, and the pressure field in each point \cite{sanderse2018energy}. The time-stepping is done with a fourth-order Runge-Kutta time-stepping scheme \cite{sanderse2012accuracy}.
The training data consists of 80 parameterized trajectories, each computed with 4000 time steps with $t\in [0,13]$. However, we train the latent time-stepping network using only every fourth time step, $s=4$, conclusively enabling the network to take longer time steps. 

The convolutional autoencoder is trained to encode the velocity field as well as the pressure field. Hence, the features are a three channel matrix input, consisting of the velocity in the $x$-direction, the $y$-direction, and the pressure field. 

For the lid driven cavity problem the same tests as for the heat equation have been made. However, the plots showing the performance for varying parameters have been omitted. To summarize, we see that the best performance is achieved with a CCNN with a memory of $\xi=8$, $\beta_1=10^{-6}$, $\beta_2=10^{-10}$, a training set of 80 trajectories, and by approximating the residual. 

\begin{figure}
\centering
\begin{subfigure}{.33\textwidth}
    \centering
    \includegraphics[width=1\textwidth]{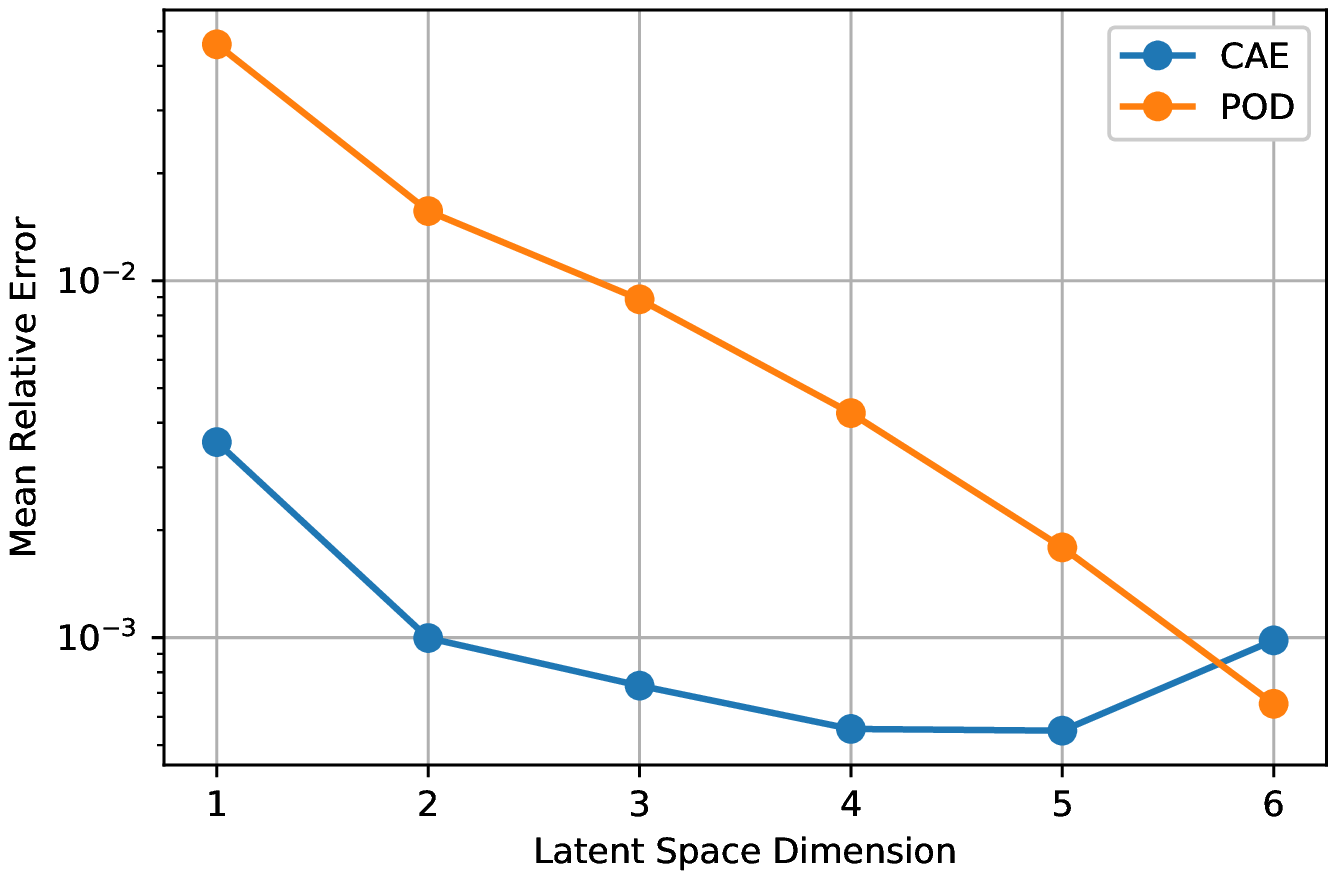}
    \caption{}
    \label{fig:pod_vs_AE_ldc}
\end{subfigure}
\begin{subfigure}{.33\textwidth}
  \centering
  \includegraphics[width=1\linewidth]{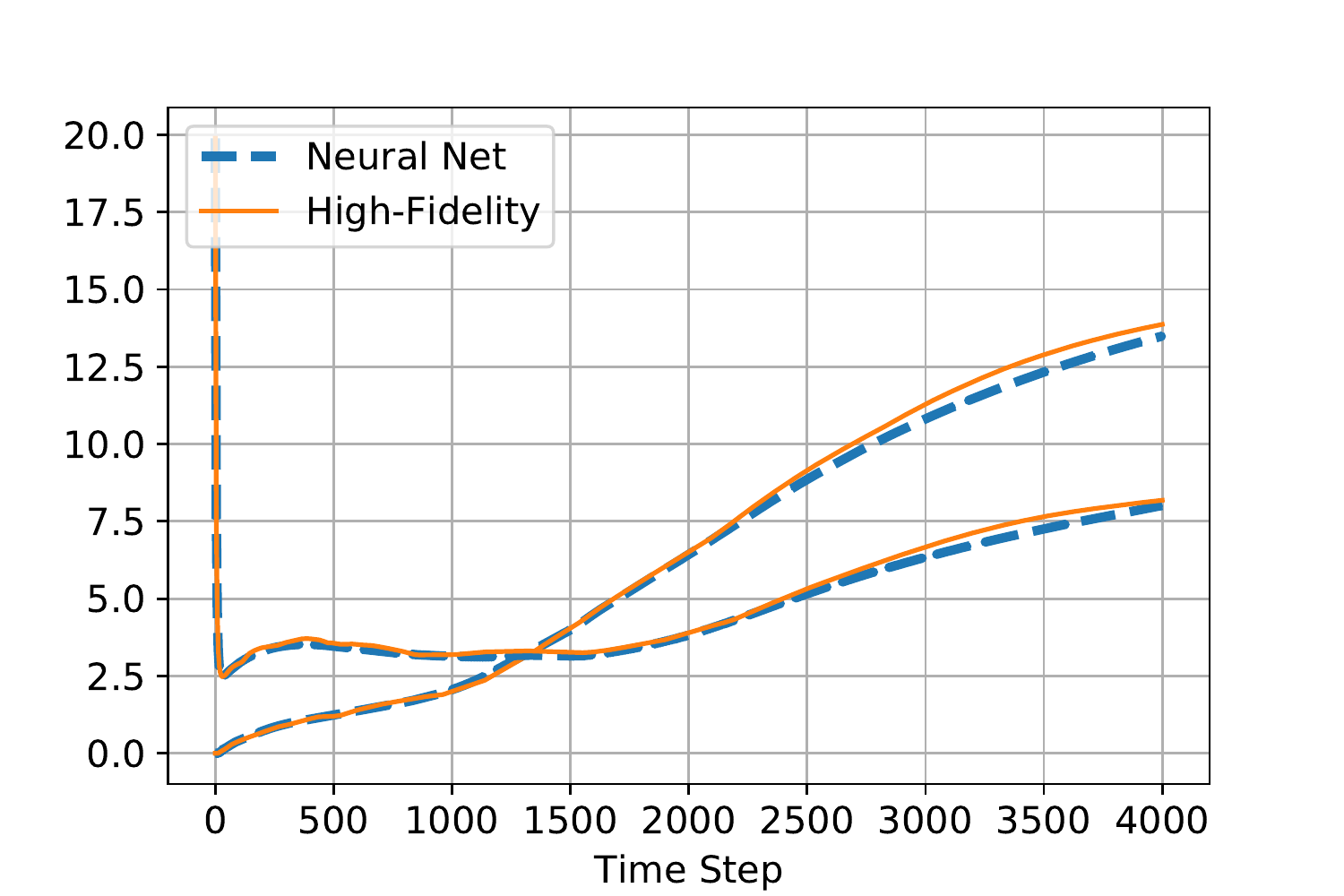}
  \caption{}
  \label{fig:reduced_ldc_sols}
\end{subfigure}%
\begin{subfigure}{.33\textwidth}
  \centering
  \includegraphics[width=1\linewidth]{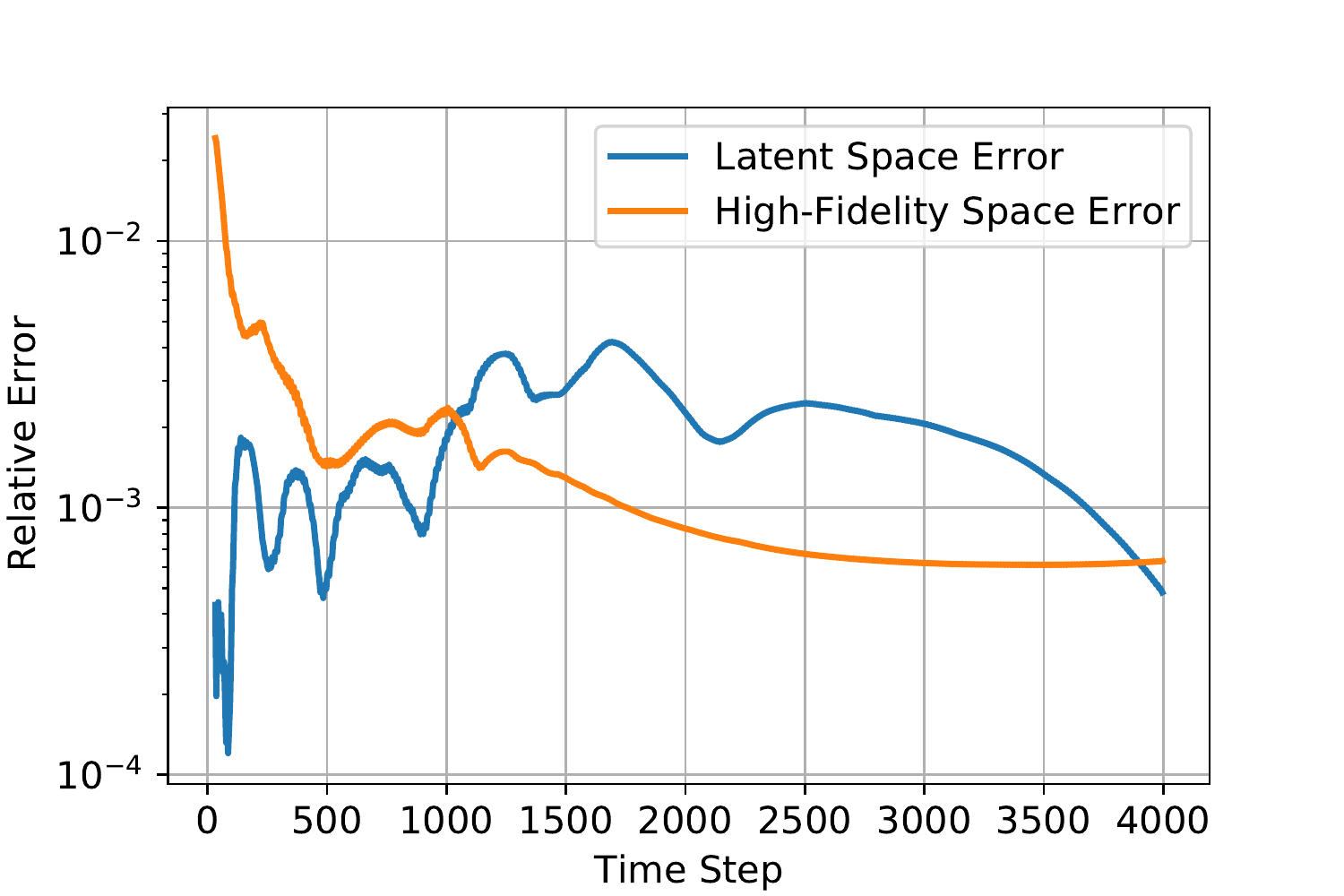}
  \caption{}
  \label{fig:ldc_time_error}
\end{subfigure}
\caption{(a) CAE and POD convergence as well as (b) trajectories and (c) errors for the lid driven cavity problem computed with Reynolds number, $Re=287$. The predictions are computed using the CCNN time-stepping network.}
\label{fig:ldc_error_time}
\end{figure}

Since we are only dealing with varying one parameter, the Reynolds number, the intrinsic dimension of the solution manifold is 2. In Figure \ref{fig:pod_vs_AE_ldc} we see that the error reaches reaches an order of magnitude at $10^{-3}$ when using a latent space dimension of exactly that. Increasing the dimension only leads to modest improvements. To achieve the same precision using POD one needs a latent space of dimension 6. Hence, we do get a significant improvement.

In Figure \ref{fig:reduced_ldc_sols} we see that the predicted latent trajectories are following the trend well. Furthermore, in Figure  \ref{fig:ldc_time_error} one sees that when the trajectories are decoded the error remains around the same order of magnitude as in the latent space, suggesting that small errors in the latent space are not inherited into the high-fidelity space.

\begin{figure}
\centering
\begin{subfigure}{.45\textwidth}
  \centering
  \includegraphics[width=.8\linewidth]{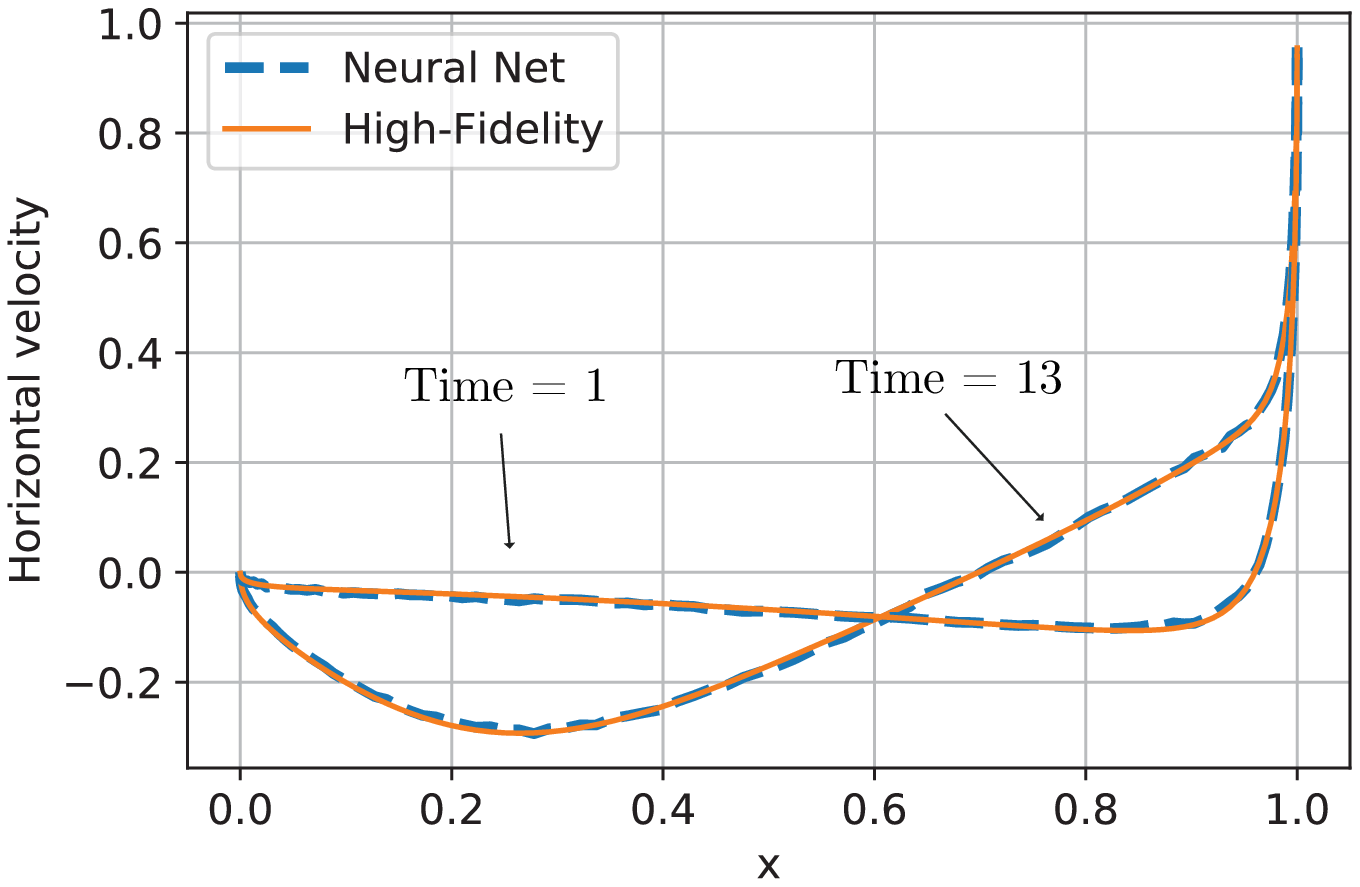}
  \caption{Horizontal velocity.}
  \label{fig:horizontal}
\end{subfigure}%
\begin{subfigure}{.45\textwidth}
  \centering
  \includegraphics[width=.8\linewidth]{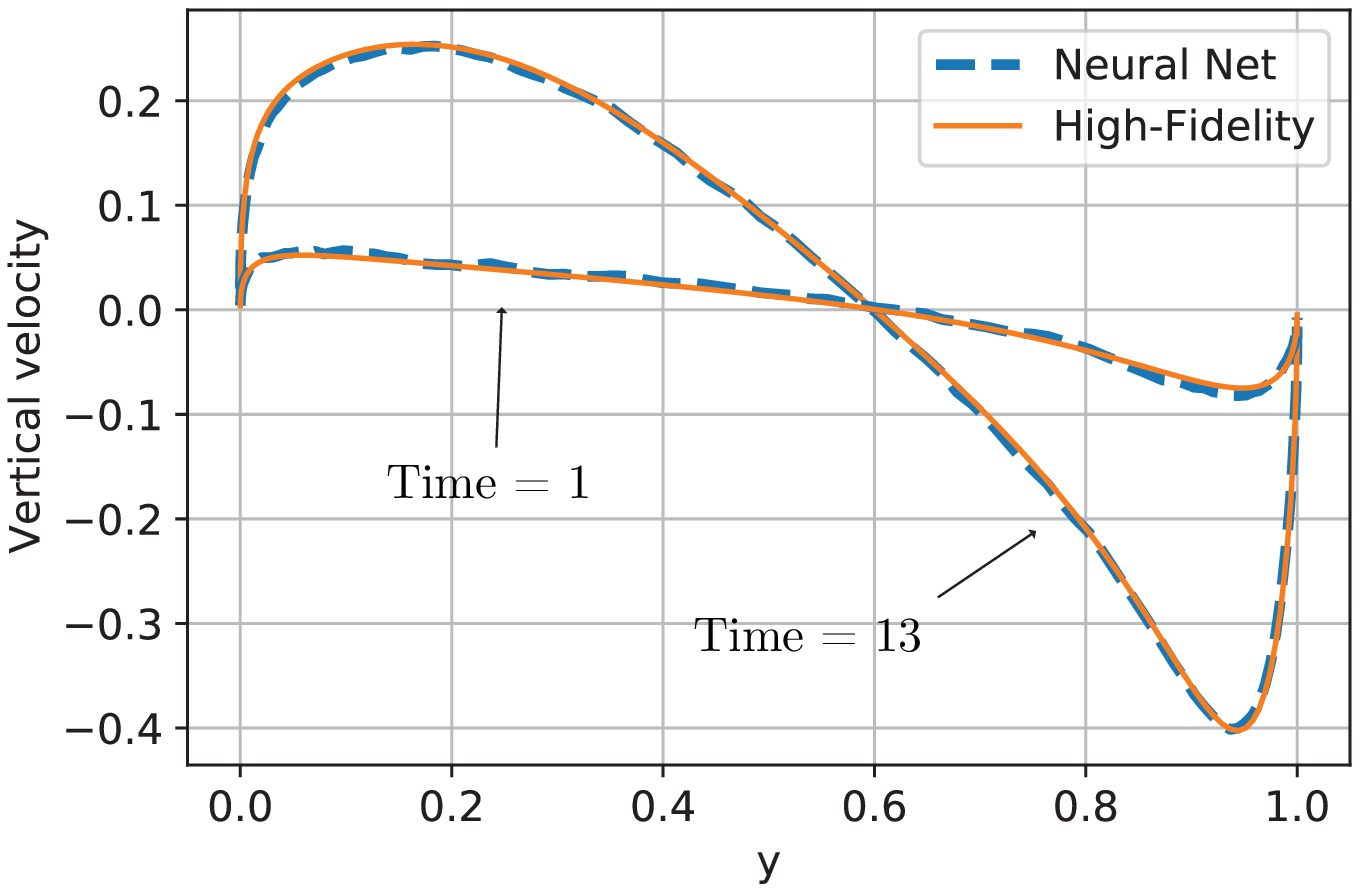}
  \caption{Vertical velocity.}
  \label{fig:vertical}
\end{subfigure}
\caption{Velocity in the horizontal and vertical direction and $0\leq x\leq1$, $y=0.5$ and $x=0.5$, $0\leq y\leq1$ respectively for $Re=287$.}
\label{fig:ldc_horizontal_vertical}
\end{figure}

When looking at common measures for the lid driven cavity, the $x$-velocity and the $y$-velocity in the along $y=0.5$ and $x=0.5$ respectively, in Figure \ref{fig:ldc_horizontal_vertical} one sees that the NN approximates the high-fidelity solution throughout the entire time interval. Even at the boundary layer we see a that the NN captures the dynamics. This is further established in Figure \ref{fig:ldc_horizontal_vertical} where one sees that the errors remain quite small during the whole timespan. 

\begin{figure}
\centering
\begin{subfigure}{.33\textwidth}
  \centering
  \includegraphics[width=1.\linewidth]{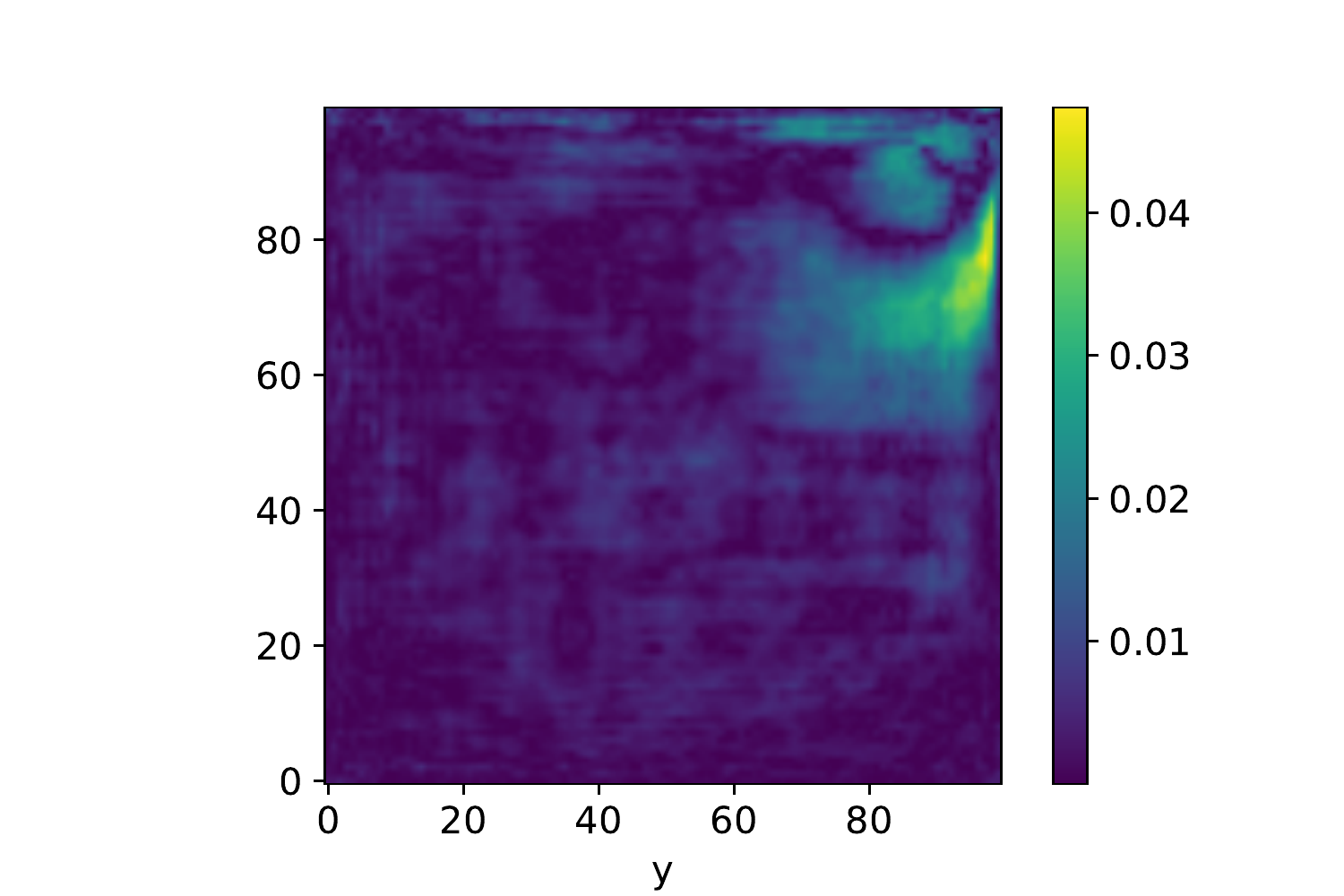}
  \caption{$t=1$}
  \label{fig:ldc_vel_magnitude_t1_error}
\end{subfigure}
\begin{subfigure}{.33\textwidth}
  \centering
  \includegraphics[width=1.\linewidth]{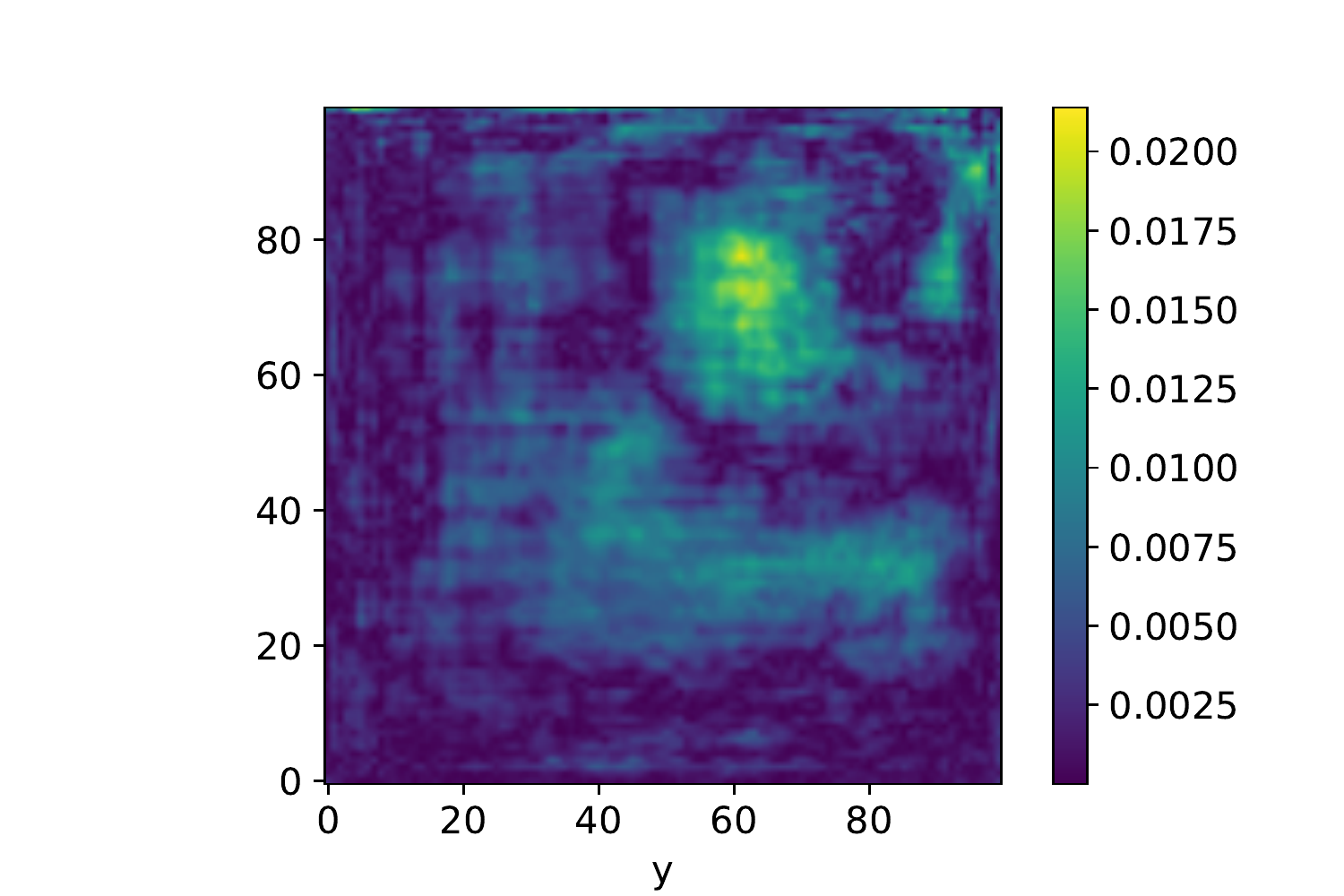}
  \caption{$t=3$}
  \label{fig:ldc_vel_magnitude_t5_error}
\end{subfigure}%
\begin{subfigure}{.33\textwidth}
  \centering
  \includegraphics[width=1.\linewidth]{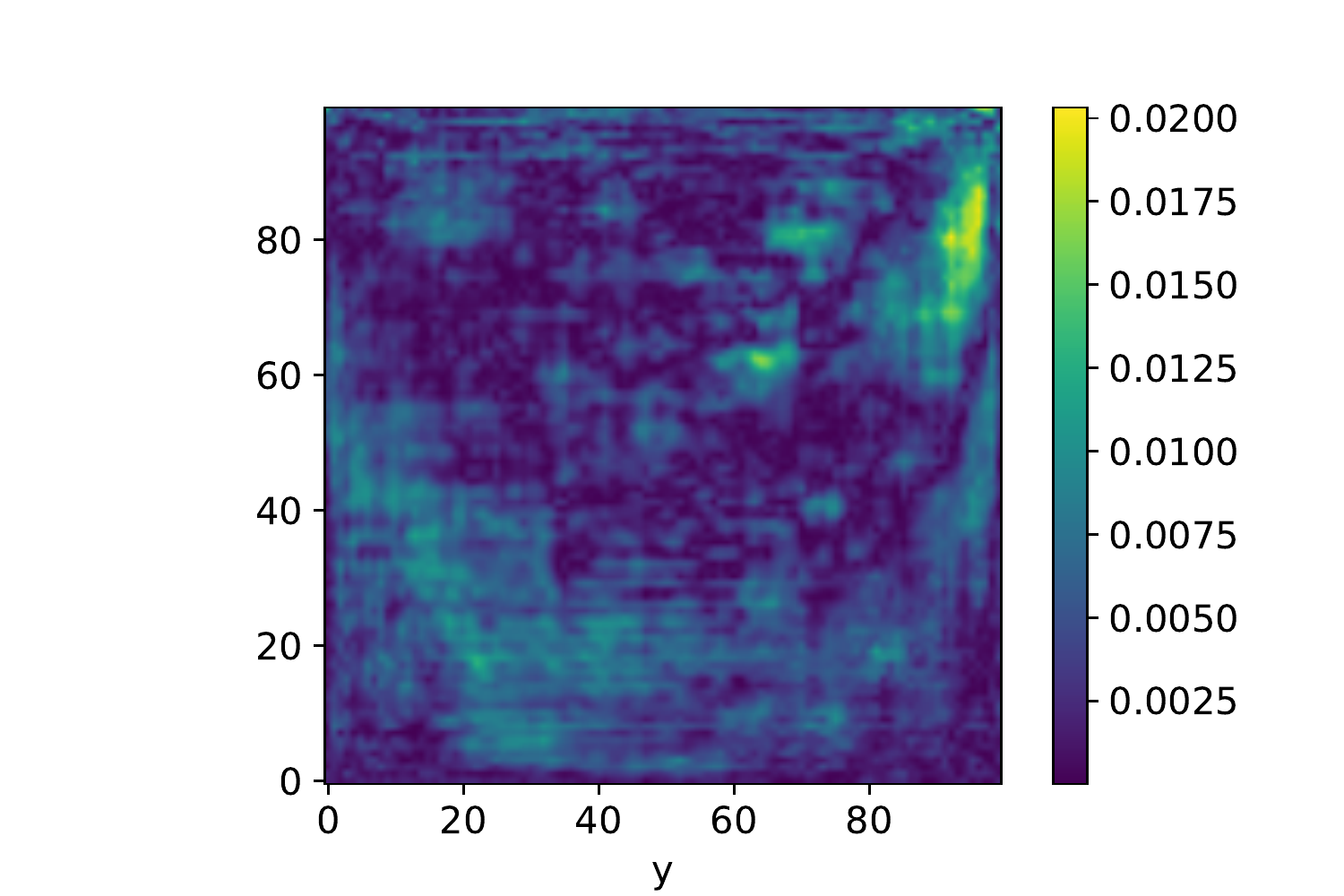}
  \caption{$t=13$}
  \label{fig:ldc_vel_magnitude_t10_error}
\end{subfigure}
\begin{subfigure}{.33\textwidth}
  \centering
  \includegraphics[width=1.\linewidth]{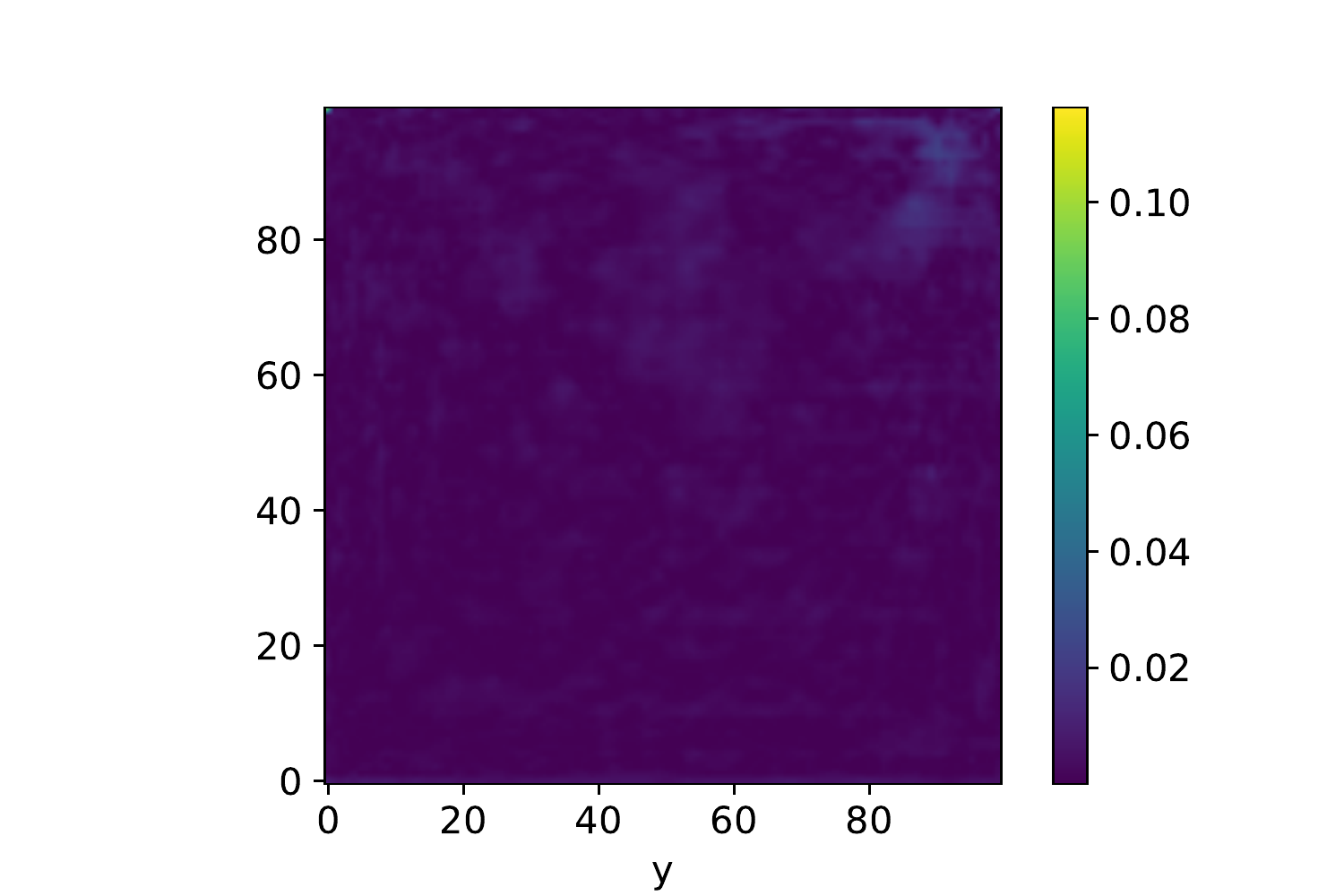}
  \caption{$t=1$}
  \label{fig:ldc_pres_magnitude_t1_error}
\end{subfigure}
\begin{subfigure}{.33\textwidth}
  \centering
  \includegraphics[width=1.\linewidth]{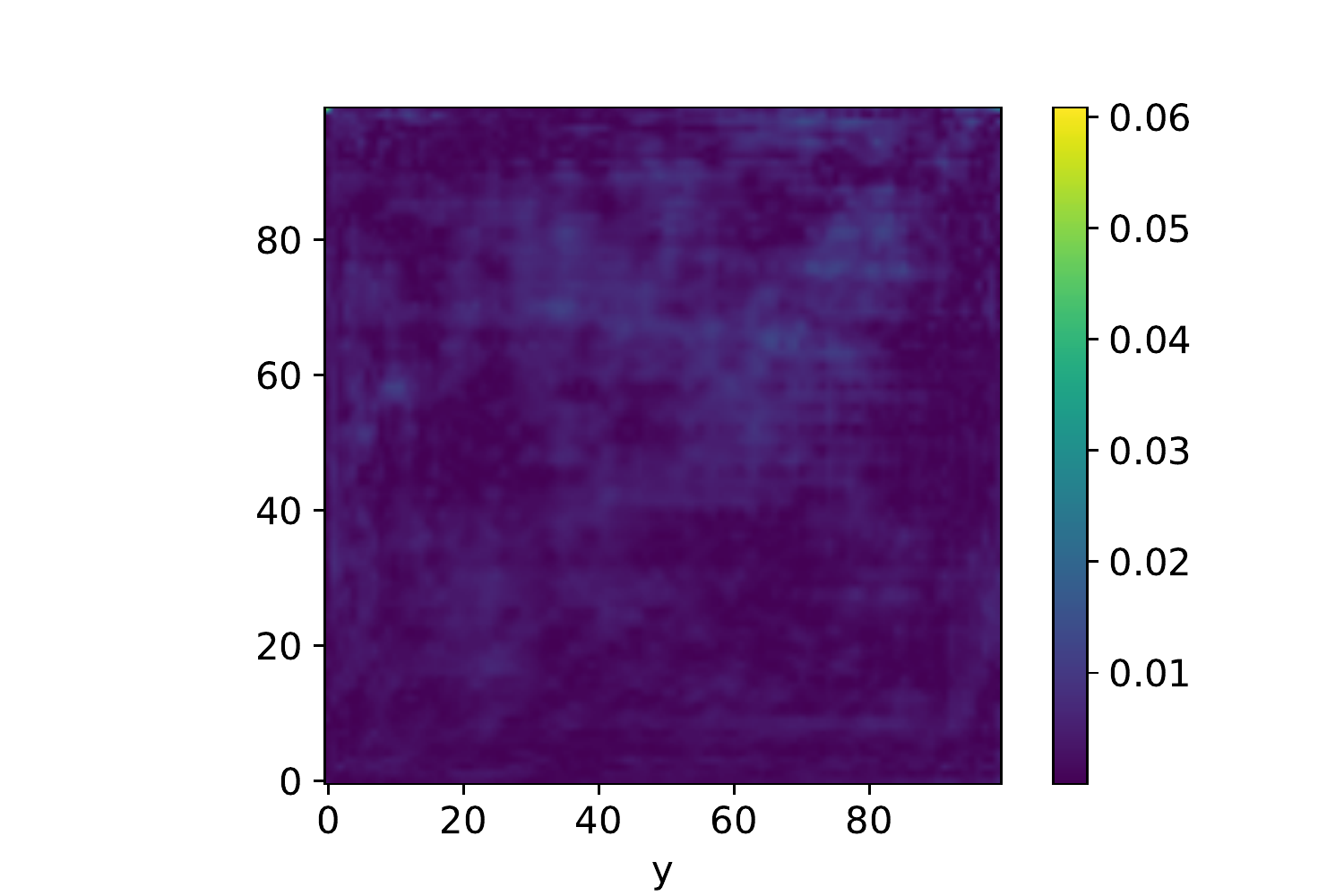}
  \caption{$t=3$}
  \label{fig:ldc_pres_magnitude_t5_error}
\end{subfigure}%
\begin{subfigure}{.33\textwidth}
  \centering
  \includegraphics[width=1.\linewidth]{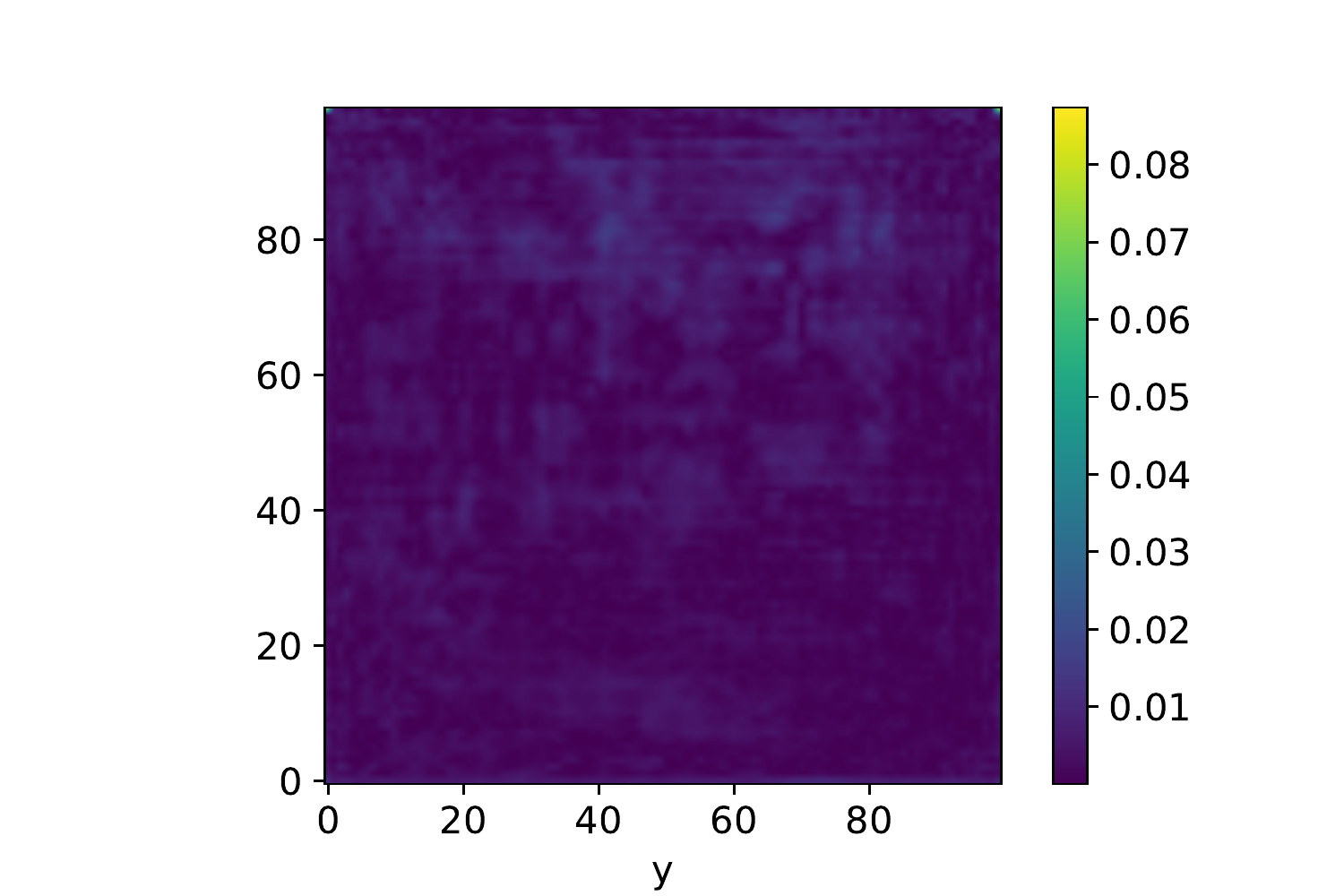}
  \caption{$t=13$}
  \label{fig:ldc_pres_magnitude_t10_error}
\end{subfigure}
\caption{Pointwise absolute error of the velocity magnitude (top row), $\sqrt{u_x^2+u_y^2}$, and the pressure field (bottom row) between the high-fidelity solution and the neural network prediction (CCNN) for the lid driven cavity problem with $Re = 287$.}
\label{fig:ldc_error}
\end{figure}

Lastly, considering the pointwise error of the velocity magnitude, Figure \ref{fig:ldc_vel_magnitude_t1_error}-\ref{fig:ldc_vel_magnitude_t10_error}, and the pointwise error of the pressure, Figure \ref{fig:ldc_pres_magnitude_t1_error}-\ref{fig:ldc_pres_magnitude_t10_error}

\subsection{Computation Time and Accuracy}
Above, we have showed and discussed performance regarding relative error for our three test cases. The results for the three cases, using the CCNN and the LSTM, are summarized in Table \ref{tab:precision}, where the time averaged error is shown. We clearly see that the differences between the LSTM and CCNN are subtle. 

As mentioned in the introduction, the aim is to be able to compute solutions fast in the online stage. In Table \ref{tab:online_time} the high-fidelity as well as the NN online time is shown. In the online stage there has not been used any form of parallelization. Therefore, it should be noted that significant speed ups for both the high-fidelity and the NN approaches could be achieved with a greater effort on this matter. The NN online time and the high-fidelity computation time is computed on an Intel Xeon 2.30GHz CPU in the CPU case and a Tesla p100 in the GPU case for the NNs. For the CCNN there is no significant difference between the online time when using a CPU compared to a GPU, while for the LSTM, the computation time approximately doubles when using a CPU.

In general, it is clear that the NN is significantly faster in the online stage. Especially, for the advection equation, we see a massive speed up. This is due to the fact that purely advection dominated problems requires relatively small time steps to avoid spurious oscillations in the solution. 

In table \ref{tab:offline_time} the offline time is shown, divided into NN training time and the time it took to generate the training trajectories. In cases where the training trajectories come from collected data the simulation step is unnecessary, and hence the training time alone is the relevant number. For the training we used a Tesla p100 GPU in Google Colab. Compared to the online stage it makes a massive difference to use a GPU instead of a CPU due to the heavy computations associated with backpropagation. We have chosen to only show the GPU training time. It is clear that the most time consuming part is generating the training trajectories. However, it should be noted that a significant amount of time has been spent on hyperparameter tuning which is not documented here. 

\begin{table}[]
    \centering
    \begin{tabular}{|c|c|c|}
    \hline
        Test Problem & CCNN & LSTM  \\ \hline
        Heat Equation & $4.41 \cdot 10^{-4}$ &  $3.55\cdot 10^{-4}$
\\
        Advection Equation & $1.54 \cdot 10^{-3}$  & $5.02\cdot 10^{-3}$
 \\
        Lid Driven Cavity & $1.56\cdot 10^{-3}$ & $1.49\cdot 10^{-3}$  \\ \hline
    \end{tabular}
    \caption{Time averaged MRE for the three test problems using the CCNN and LSTM.}
    \label{tab:precision}
\end{table}

\begin{table}[]
\centering
\begin{tabular}{ c c||c|c||c|c| }
 \cline{2-6}
 \multicolumn{1}{c||}{} & \multicolumn{1}{c||}{High-Fidelity} &  \multicolumn{2}{c||}{CCNN}  & \multicolumn{2}{c|}{LSTM} \\
 \hline
 \multicolumn{1}{|c||}{Test Problem}&\multicolumn{1}{c||}{CPU} & CPU & GPU & CPU & GPU\\ 
 \hline
 \multicolumn{1}{|c||}{Heat Equation}   & 42.62  & 0.58  & 0.48 & 1.21 & 0.49    \\
 \multicolumn{1}{|c||}{Linear Advection Equation} & 532.22 & 2.08 & 2.03 & 4.63  & 2.39\\
 \multicolumn{1}{|c||}{Lid Driven Cavity} & 94.83  & 5.15 & 4.54 & 11.18 & 4.52 \\
 \hline
\end{tabular}
 \caption{Online computationn. We remind the reader that the following schemes are used for the high-fidelity computations: i) Third-order Galerkin finite element method on a $100\times 100$ grid and a Crank-Nicolson time-stepping scheme for the heat equation, ii) second-order discontinuous Galerkin finite element method on a $60\times 60$ grid and a Crank-Nicolson time-stepping scheme for the advection equation, and iii) a second-order finite volume scheme on a $100\times 100$ non-uniform staggered grid and a fourth-order explicit Runge-Kutta time-stepping scheme for the lid driven cavity problem. }
    \label{tab:online_time}
\end{table}

\begin{table}[]
    \centering
    \begin{tabular}{|c|c|c|c|} \hline
        Test Problem &  CCNN & LSTM & Generation of Trajectories \\ \hline
        Heat Equation &  1178.89 & 1273.63  & 11485.50 (1450 trajectories) \\
        Advection Equation &  396.67 & 622.28  & 50763.13 (120 trajectories) \\
        Lid Driven Cavity & 798.00 & 761.59 & 7121.4 (80 trajectories) \\ \hline
    \end{tabular}
    \caption{Offline computation time, i.e. NN training time, in seconds for the CCNN and LSTM using GPUs. Furthermore we show the time it took to generate the training trajectories. Note that the generation of training trajectories is not necessary in cases where the data already exists. }
    \label{tab:offline_time}
\end{table}

%% file: Conclusion.tex
\section{Conclusion} \label{Conclusion}

We presented a novel deep learning approach to non-intrusive reduced order modeling for parameterized time-dependent PDEs using CAEs for dimensionality reduction and CCNNs and LSTMs combined with FFNNs for time evolution. This approach was demonstrated on various test cases and was shown to perform well in the online phase, showcasing the potential of using deep learning based ROMs for different physical phenomena. 

Regarding dimensionality reduction, a discussion and comparison of linear and nonlinear methods was presented with POD and CAEs as the focus points. The discussion focused on why a nonlinear approach has the potential to outperform a linear approach. 

For time stepping, the general idea was to encode the previous states and the parameters separately in parallel and then combine the encoded data to make a final prediction using a FFNN. The two encoding NNs, as well as the final prediction NN, constitute a single network, meaning everything is trained simultaneously. This ensures that both the memory and parameters are encoded in relation to one another. Furthermore, various methods to ensure generalization, stability, and precision were discussed and tested. 

In all the test cases errors were found to be between $10^{-4}$ and $10^{-2}$ for all time steps. These are in many cases acceptable errors considering the significant speed-ups. Furthermore, we saw that the model performs well on both rather high-dimensional trial manifolds (the heat equation), purely advection dominated problems (the advection equation) and nonlinear multiple vector field computations (the lid driven cavity) and it delivered significant online speed-ups without sacrificing accuracy. 
 
In summary, the contributions in this work include a nonlinear dimensionality reduction scheme using convolutional autoencoders, a novel parallel neural network architecture for parameterized time-stepping using CCNNs and LSTMs, and a discussion on different approaches to achieve stability and generalization for neural network-based time-stepping.

In the future the methodology will be tested on more advanced PDE problems. By advanced problems, we are both referring to increasing nonlinearity, higher dimensions, and multi-query problems such as uncertainty quantification, model predictive control, and data assimilation. 

Besides considering other use cases one could work on improving the NN architecture and training by, e.g. incorporating the physics in the training \cite{erichson2019physics, raissi2017physics}, and use Bayesian optimization \cite{archetti2019bayesian} or reinforcement learning \cite{sugiyama2015statistical} to ensure effective snapshot generation. Furthermore, with the amount of hyperparameters ($\xi$, $\beta_1$, $\beta_2$, number of layers and neurons, etc.) the task of hyperparameter tuning is not trivial and could be solved more effectively using modern approaches. 

\section*{Acknowledgement}
This work is supported by the Dutch National Science Foundation NWO under the grant number 629.002.213.
which is a cooperative projects with IISC Bangalore and Shell research as project partners.
The authors furthermore acknowledge fruitful discussions with Dr. B. Sanderse.

%% file: Artificial_Neural_Networks.tex
\section{Artificial Neural Networks} \label{appendix:neural_network}

\subsection{Feedforward Neural Networks}
The arguably most common ANN architecture is the feedforward neural network (FNN). An FNN can be considered a function, $G:\R^{N_i}\rightarrow \R^{N_o}$, consisting of a series of affine transformations, $T_i$, followed by an element-wise (nonlinear) activation function, $\sigma_i$:
\begin{align} \label{FNN}
    G(x;\theta) = \sigma_L\circ T_L\circ \ldots \circ \sigma_1 \circ T_1 (x).
\end{align}
The combination of an afine transformation followed by the activation is called a neuron. The afine transformation can be written as $T_i(x)=W_ix+b_i$, where $W\in \R^{M_i\times M_{i-1}}$ and $b\in\R^{M_i}$. We call $W_i$ the weight matrix, $b_i$ the bias vector, and $M_i$ the number of neurons in layer $i$, and $L$ the number of layers. \eqref{FNN} is conveniently visualized as a network of neurons. We will refer to the set of parameters as $\theta = \left\{W_1,b_1, \ldots,W_L,b_L\right\}$.

\begin{figure}
    \centering
    \includegraphics[width=.35\textwidth]{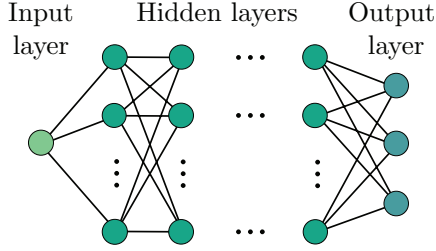}
    \caption{Visualization of a feedforward densely connected neural network.}
    \label{fig:NN_dense}
\end{figure}

In supervised learning one tries to approximate a function by an ANN, typically done by minimizing the empirical risk w.r.t. the parameters $\theta$:
\begin{align} \label{minimization}
    \theta^* = \text{arg}\min_\theta \mathbb{E}_{(x,y)\sim \mathcal{P}_{\text{data}}} \left[ \mathcal{L}(G(x;\theta),y)  \right],
\end{align}
where $\mathcal{P}_{\text{data}}$ is the distribution generating the data and $\mathcal{L}$ is a chosen loss function measuring the discrepancy between the predicted output and the target. For regression type problems the mean squared error (MSE) is the most common choice. However, especially for physics-informed machine learning the physics is often incorporated in the loss function in shape of extra term \cite{raissi2017physics, erichson2019physics}.

Computing \eqref{minimization}, i.e. training the ANN, is mostly done using stochastic gradient descent (SGD) or a variants such as ADAM \cite{kingma2014adam}.

\subsection{Convolutional Neural Networks}
Convolutional neural networks (CNNs) gained attention due their great performance in image recognition. The general idea is to utilize local properties of the data instead of only considering global properties. This is done by having local connections and shared weights in the neural networks. These properties are not only great for detecting patterns in data but it also makes it possible to do computations on very high-dimensional data. 

A convolutional layer is effectively a feature map where each unit in the layer is connected to a local patch of the previous layer through a filter bank and an activation function. A feature map at layer $l$ is a tensor, $H^l\in\R^{N_{chan}^l\times N_1^l \times N_2^l}$, where $H^l_{i,j,k}$ is a unit at channel $i$, row $j$, and column $k$. The filter bank at layer $l$ is a 4-dimensional tensor, $F^l\in \R^{N_{filter}^l\times N_{chan}^{l-1}\times k_1 \times k_2}$, where $F^l_{i,j,m,n}$ connects a unit in channel $i$ of the output and channel $j$ of the input with $m$ and $n$ being the offset of rows and columns respectively. $N_{filter}^l$ denotes the number of filters in the feature bank in layer $l$ and $k_1$ and $k_2$ denotes the kernel size. The convolution operation between a feature map and a filter bank is given by
\begin{align}
    H^l_{i,j,k} = \sigma_l\left(\sum_{r=1}^{N_{chan}^{l-1}}\sum_{m=1}^{N_1^{l-1}}\sum_{n=1}^{N_2^{l-1}}
    H^{l-1}_{r,(j-1)s+m,(k-1)s+n} F_{i,r,m,k}^{l} + B^l_{i,j,k} \right),
\end{align}
where $B^l_{i,j,k}$ is a bias term and $\sigma_l$ is an activation function applied element-wise. $s$ denotes the stride and effectively downsamples the feature map between layers. The filters, $F_{i,r,m,k}^{l}$, and biases, $B^l_{i,j,k}$, are the learnable parameters while the kernel sizes, $k_1, k_2$, the stride, $s$, and the number of filters, $N_{filter}^l$, are chosen. Often these are subject case specific objectives or hyperparameter optimization. 

\subsection{Causal Convolutional Neural Networks}
As the name suggests, causal convolutional neural networks (CCNNS) are related to convolutional neural networks. CCNNs are sometimes referred to as temporal convolutional neural networks, but in this paper we use the term CCNN. 

CCNNs are used for encoding time series data with the purpose of forecasting or classification. The general idea is to use 1-dimensional convolutions on time series data. In the multivariate case the multiple dimensions are interpreted as channels. The term causal refers to the fact that the filter banks are only convolved with the current and previous time steps, thus establishing a causal relationship between the past the future. 

\subsection{Recurrent Neural Networks and Long Short-Term Memory}
A recurrent neural network (RNN) is an alternative to CCNNs for interpreting time series data. The general idea is to process and retain information from previous time steps in an efficient manner. In this paper, we solely focus on a specific RNN called long short-term memory (LSTM) \cite{hochreiter1997long}. For an input consisting of several previous time steps, $x^n$, an LSTM layer consists of four components \cite{gonzalez2018deep}: An input gate:
\begin{align}
    i^{n+1} = \sigma\left(W_i x^{n} + b_i \right),
\end{align}
a forget gate:
\begin{align}
    f^{n+1} = \sigma\left(W_f x^{n} + b_f \right),
\end{align}
an output gate:
\begin{align}
    o^{n+1} = \sigma\left(W_o x^{n} + b_o \right),
\end{align}
and a cell state
\begin{align}
    c^{n+1} = i \odot c^{n} + i^n \odot \tanh\left(W_c x^{n} + b_c \right).
\end{align}
The prediction is then given by
\begin{align}
    x^{n+1} = o^n \odot \tanh\left(c^n\right).
\end{align}
$W_i,b_i,W_f,b_f,W_o,b_o,W_c,b_c$ are the trainable weight matrices and bias vectors, and $\odot$ is the Hadamard product. Ideally, the input gate identifies what information to be passed to the from the cell state, the forget what to be dropped, and the output gate decides what to be passed to the final prediction.

%% file: appendix_results.tex
\section{Convolutional Autoencoder Configurations} \label{conv_ae_appendix}

\begin{table}[H]
    \centering
    \begin{tabular}{c|c|c|c}
        \textbf{Encoder} & & & \\ \hline \hline
        \textbf{Convolutional Layers} & Kernels & Filter Size & Stride \\ \hline
        Convolutional Layer 1 & 4 & $5\times 5$ & $2\times 2$   \\
        Batchnormalization 1 & & &  \\
        Convolutional Layer 2 & 8 & $5\times 5$ & $2\times 2$   \\
        Batchnormalization 2 & & &  \\
        Convolutional Layer 3 & 16 & $5\times 5$ & $2\times 2$   \\
        Batchnormalization 3 & & &  \\
        Convolutional Layer 4 & 32 & $5\times 5$ & $2\times 2$   \\
        Batchnormalization 4 & & &  \\
        Flatten & & & \\ \hline
        \textbf{Dense Layers} & Neurons & & \\ \hline 
        Dense Layer 1 & $N_l$ & & \\ \hline
        & & & \\ \hline
        \textbf{Decoder} & & & \\ \hline \hline
        \textbf{Dense Layers} & Neurons & Size & \\ \hline 
        Dense Layer 1 & 512 & & \\
        Reshape & & $4\times 4\times 32$ & \\ \hline
        \textbf{Convolutional Layers} & Kernels & Filter Size & Stride \\ \hline
        Transposed Convolutional Layer 1 & 16 & $5\times 5$ & $2\times 2$   \\
        Batchnormalization 1 & & &  \\
        Transposed Convolutional Layer 2 & 8 & $5\times 5$ & $2\times 2$   \\
        Batchnormalization 2 & & &  \\
        Transposed Convolutional Layer 3 & 4 & $5\times 5$ & $2\times 2$   \\
        Batchnormalization 3 & & &  \\
        Transposed Convolutional Layer 5 & 1 & $5\times 5$ & $2\times 2$  \\ \hline
    \end{tabular}
    \caption{Convolutional autoencoder configuration for the heat equation anf advection equation.}
    \label{tab:con_ae_configuration_heat_adv}
\end{table}

\begin{table}[H]
    \centering
    \begin{tabular}{c|c|c|c}
        \textbf{Encoder} & & & \\ \hline \hline
        \textbf{Convolutional Layers} & Kernels & Filter Size & Stride \\ \hline
        Convolutional Layer 1 & 4 & $5\times 5$ & $2\times 2$   \\
        Batchnormalization 1 & & &  \\
        Convolutional Layer 2 & 8 & $5\times 5$ & $2\times 2$   \\
        Batchnormalization 2 & & &  \\
        Convolutional Layer 3 & 16 & $5\times 5$ & $2\times 2$   \\
        Batchnormalization 3 & & &  \\
        Convolutional Layer 4 & 32 & $5\times 5$ & $2\times 2$   \\
        Batchnormalization 4 & & &  \\
        Flatten & & & \\ \hline
        \textbf{Dense Layers} & Neurons & & \\ \hline 
        Dense Layer 1 & $64$ & & \\ 
        Dense Layer 2 & $32$ & & \\ 
        Dense Layer 3 & $N_l$ & & \\ \hline
        & & & \\ \hline
        \textbf{Decoder} & & & \\ \hline \hline
        \textbf{Dense Layers} & Neurons & Size & \\ \hline 
        Dense Layer 1 & 32 & & \\
        Dense Layer 2 & 64 & & \\
        Dense Layer 3 & 1568 & & \\
        Reshape & & $7\times 7\times 32$ & \\ \hline
        \textbf{Convolutional Layers} & Kernels & Filter Size & Stride \\ \hline
        Transposed Convolutional Layer 1 & 16 & $5\times 5$ & $2\times 2$   \\
        Batchnormalization 1 & & &  \\
        Transposed Convolutional Layer 2 & 8 & $5\times 5$ & $2\times 2$   \\
        Batchnormalization 2 & & &  \\
        Transposed Convolutional Layer 3 & 4 & $5\times 5$ & $2\times 2$   \\
        Batchnormalization 3 & & &  \\
        Transposed Convolutional Layer 5 & 2 & $5\times 5$ & $2\times 2$  \\ \hline
    \end{tabular}
    \caption{Convolutional autoencoder configuration for the lid driven cavity problem.}
    \label{tab:con_ae_configuration_ldc}
\end{table}

\section{Time Evolution Neural Network Configuration} \label{evolution_NN_appendix}

For the CCNN memory encoding the layers are organized as shown in Figure \ref{fig:CCNN_types}. For the LSTM we work network architectures of 3 layers 32 neurons in each LSTM layer. 

Furthermore, before the LSTM or CCNN layers a each previous state is passed through a dense layer with 16 neurons. In TensorFlow 2.0 this type of layer is denoted \texttt{TimeDistributed}. 

For the parameter encoding the neural network is a simple 3 layer deep network with 16 neurons in each layer.

For the final prediction we utilize a 3 layer deep NN with 32 neurons in each layer.